\documentstyle [openbib,fullpage,12pt]{amsart}

\newcommand{\qaa}{{q ; a_1, \dots , a_r}}
\newcommand{\qnn}{{8; 3, 5, 7}}
\newcommand{\sumprime}{\mathop{\sum\nolimits'}}
\newcommand{\R}{{\bf R}}
\newcommand{\Z}{{\bf Z}}
\newcommand{\PV}{{\rm P.V.}}
\newcommand{\li}{\mathop{\rm li}}

\renewcommand{\Re}{\mathop{\rm Re}}
\renewcommand{\Im}{\mathop{\rm Im}}

\newtheorem{theorem}{Theorem}
\newtheorem{lemma}{Lemma}[section]

\begin{document}

\title{Biases in the Shanks--R\'enyi Prime Numbers Race}
\author{Andrey Feuerverger\, and\, Greg Martin\\
University of Toronto}
\maketitle

{\smaller\narrower

\centerline{\sc Abstract}\medskip

Rubinstein and Sarnak investigated systems of inequalities of the form
$\pi(x;q,a_1)>\dots>\pi(x;q,a_r)$, where $\pi(x;q,a)$ denotes the
number of primes up to $x$ that are congruent to $a \bmod q$. They
showed, under standard hypotheses on the zeros of Dirichlet
$L$-functions mod $q$, that the set of positive real numbers $x$ for
which these inequalities hold has positive (logarithmic) density
$\delta_{q;a_1,\dots,a_r}>0$. They also discovered the surprising fact
that a certain distribution associated with these densities is not
symmetric under permutations of the residue classes $a_i$ in general,
even if the $a_i$ are all squares or all nonsquares mod $q$ (a
condition necessary to avoid obvious biases of the type first observed
by Chebyshev). This asymmetry suggests, contrary to prior
expectations, that the densities $\delta_{q;a_1,\dots,a_r}$ themselves
vary under permutations of the~$a_i$.

In this paper, we derive (under the hypotheses used by Rubinstein and
Sarnak) a general formula for the densities
$\delta_{q;a_1,\dots,a_r}$, and we use this formula to calculate many
of these densities when $q\le12$ and $r\le4$. For the special moduli
$q=8$ and $q=12$, and for $\{a_1,a_2,a_3\}$ a permutation of the
nonsquares $\{3,5,7\}\bmod8$ and $\{5,7,11\}\bmod12$, respectively, we
rigorously bound the error in our calculations, thus verifying
that these densities are indeed asymmetric under permutation of the
$a_i$. We also determine several situations in which the densities
$\delta_{q;a_1,\dots,a_r}$ remain unchanged under certain permutations
of the $a_i$, and some situations in which they are provably
different.

\medskip

\noindent {\it Key words and phrases:} Chebyshev's bias, comparative
prime number theory, primes in arithmetic progressions, Shanks-R\'enyi
race.

}

\medskip
\section{Introduction and Summary.}

In 1853 Chebyshev remarked that there are more primes congruent to 3
than to 1 modulo~4, and since that time considerable efforts have been
expended in attempts to determine in what sense this remark is true.
It follows from the prime number theorem for arithmetic progressions
(see for instance Davenport \cite{davenport}) that, asymptotically,
half of all primes are congruent to $3 \bmod 4$ and half are congruent
to $1 \bmod 4$, so that Chebyshev's observation cannot be interpreted
in that sense.  However, when we compute the numbers of primes up to
$x$ that are congruent to $3 \bmod 4$ and to $1 \bmod 4$, we find that
for most values of $x$, the primes congruent to 3 are more numerous
than those congruent to 1.  Similar ``biases'' have also been
observed, notably by Shanks~\cite{shanks}, for moduli $q$ other than
4; in particular, the numbers of primes in nonsquare residue classes
modulo $q$ tend to exceed the numbers of primes in square residue
classes. We refer to inequities of this type as ``Chebyshev biases''.

These observations lead naturally to the study of inequalities of the type
\begin{equation}
\pi(x;q,a_1)  >  \pi(x;q,a_2)  >  \dots  >  \pi(x;q,a_r),
\label{pi.inequality}
\end{equation}
where $\pi(x;q,a)$ denotes the number of primes $p \leq x$ such that $
p \equiv a \bmod q$.  Littlewood~\cite{littlewood} showed
(unconditionally) that the inequalities $\pi(x;3,1) > \pi(x;3,2)$ and
$\pi(x;4,1) > \pi(x;4,3)$, as well as the opposite inequalities, each hold
for infinitely many integer values of $x$.  A number of additional
results on single inequalities of this type were subsequently derived
under certain hypotheses by Knapowski and Tur\'an in a series of
papers beginning with~\cite{KT}, and Kaczorowski wrote several
papers concerning the multiple inequalities~(\ref{pi.inequality}),
the most recent of which is~\cite{kaczorowski}.

A major advance was made recently by Rubinstein and Sarnak \cite{RS}
who showed (conditionally) that for any modulus $q$ and for any distinct
reduced residues $a_1, \dots , a_r \bmod q$ (i.e., integers relatively
prime to $q$), the system of inequalities (\ref{pi.inequality}) holds
for infinitely many integers $x$.  More precisely, they worked under
the assumption of the Generalized Riemann Hypothesis for Dirichlet
$L$-functions, which we shall abbreviate GRH, and an additional
assumption (their ``Grand Simplicity Hypothesis'') that the imaginary
parts of the nontrivial zeros of Dirichlet $L$-functions corresponding
to primitive characters are linearly independent over the rationals,
which we shall abbreviate LI. Rubinstein and Sarnak studied the
quantities $\delta_\qaa$, defined as the logarithmic density of the
set of positive real numbers $x$ for which the system of inequalities
(\ref{pi.inequality}) holds. (Here, the logarithmic density $\delta
(\Lambda)$ of any subset $\Lambda$ of the real numbers is defined as
\begin{equation*}
\delta (\Lambda) \  = \  \lim_{x \rightarrow \infty} \  \frac{1}{\log x}
\int_{\Lambda \cap [2,x]} \frac{dt}{t} \, ,
\end{equation*}
provided that this limit exists. Suffice it to say here that
logarithmic densities are more appropriate for these problems than
ordinary densities; in this paper, by ``density'' we shall always mean
logarithmic density.)

Under the above hypotheses, Rubinstein and Sarnak proved that the
densities $\delta_\qaa$ exist and are positive for any integer $q\ge2$
and for any distinct reduced residues $a_1, \dots , a_r$ mod $q$.
They obtained, for several small moduli $q$, numerical values for the
density of those $x$ for which
the primes up to $x$ that are quadratic nonresidues mod $q$
outnumber those which are quadratic residues. Rubinstein and Sarnak
also proved that $\delta_{q;a,a'} = \delta_{q;a',a} = 1/2$ if $a$ and
$a'$ are both squares or both nonsquares mod $q$, and otherwise
$\delta_{q;a,a'}$ is greater than or less than $1/2$ according to
whether $a$ or $a'$ is the nonsquare mod $q$, thus bearing out the
biases of the type observed by Chebyshev.

It was generally suspected for $r>2$ as well that whenever the $a_j$
are all squares or all nonsquares modulo $q$, the densities
$\delta_\qaa$ are invariant under permutations of the $a_j$ (and thus
equal to $1 / r!$). However, Rubinstein and Sarnak showed that certain
distributions $\mu_\qaa$ on $\R^r$ that are associated naturally with
the densities $\delta_\qaa$ are not symmetric under permutations of
the $a_j$ when $r \geq 3$, except in the special case when $r=3$ and
there exists $\rho \not\equiv 1 \bmod q$ with $\rho^3 \equiv 1 \bmod
q$ such that $a_2 \equiv {a_1\rho} \bmod q$ and $a_3 \equiv
{a_1\rho^2} \bmod q$. (Note that since $\rho\equiv\rho^4\bmod q$ is a
square, it follows that such $\{a_1, a_2, a_3\}$ are all squares or
all nonsquares mod $q$.) This result suggests, but does not imply,
that the $\delta_\qaa$ are generally asymmetric under permutation of
the $a_j$.

In this paper, we rigorously establish a number of asymmetries of this
type.  Triples of nonsquares and triples of squares occur for the
moduli $q=7$ and $q=9$, but these triples fall under the special case
that has just been mentioned.  Therefore the smallest moduli for which
such asymmetries of the $\delta_\qaa$ could arise are $q=8$ and
$q=12$, each of which has three nonsquares (and a single square),
and $q=11$, which has five squares and five nonsquares.  Our main theorem
provides rigorous results for the cases $q=8$ and $q=12$, subject to
the two aforementioned hypotheses:

\begin{theorem}
Assume GRH and LI. Let $\delta_\qaa$ denote the
(logarithmic) density of the set of positive real numbers $x$ for
which the system of inequalities (\ref{pi.inequality}) holds. Then
\begin{equation*}
\begin{split}
\delta_\qnn &= \delta_{8;7,5,3} =  0.1928013 \pm 0.000001
\\
\delta_{8;3,7,5} &= \delta_{8;5,7,3} =  0.1664263 \pm 0.000001
\\
\delta_{8;5,3,7} &= \delta_{8;7,3,5} =  0.1407724 \pm 0.000001
\end{split}
\end{equation*}
and
\begin{equation*}
\begin{split}
\delta_{12;5,7,11} &= \delta_{12;11,7,5} =  0.1984521 \pm 0.000001
\\
\delta_{12;5,11,7} &= \delta_{12;7,11,5} =  0.1215630 \pm 0.000001
\\
\delta_{12;7,5,11} &= \delta_{12;11,5,7} =  0.1799849 \pm 0.000001,
\end{split}
\end{equation*}
where the indicated error bounds are rigorous.
\end{theorem}

The pairwise equalities among the $\delta$'s in Theorem 1 are not
numerical coincidences, but are provably exact. In fact there are
several situations in which we can establish symmetries of this
sort. To state these results, we first need to define
\begin{equation}
c(q,a) = {-1} + \# \{ 1 \leq b \leq q\colon b^2 \equiv a \bmod q \}
\label{cqadef}
\end{equation}
for coprime integers $a$ and $q$. Note that when $q$ is an odd prime,
$c(q,a)$ simply equals the Legendre symbol $\big( \frac{a}{q}
\big)$. Note further that $c(q,a)$ can take only two possible values
for a given~$q$: certainly $c(q,a)={-1}$ for every nonsquare $a$ mod
$q$, while $c(q,a)=c(q,1)$ for every square $a$ mod $q$. We can
interpret $c(q,1)$ as the ratio of the number of invertible nonsquares
to the number of invertible squares mod $q$.

We may now state our results concerning symmetries:

\renewcommand{\labelenumi}{(\alph{enumi})}
\begin{theorem}
Assume GRH and LI. Let $q,r\ge2$ be integers and let
$a_1,\dots,a_r$ be distinct reduced residue classes mod $q$.

\begin{enumerate}
\item Letting $a_j^{-1}$ denote the multiplicative inverse of $a_j$
modulo $q$, we have $\delta_{q;a_1,\dots,a_r} =
\delta_{q;a_1^{-1},\dots,a_r^{-1}}$.

\item If $b$ is a reduced residue class modulo $q$ such that
$c(q,a_j)=c(q,ba_j)$ for each $1\le j\le r$, then
$\delta_{q;a_1,\dots,a_r} = \delta_{q;ba_1,\dots,ba_r}$. In
particular, this holds if $b$ is a square modulo $q$.

\item If the $a_j$ are all squares modulo $q$ and $b$ is any reduced
residue class modulo $q$, then $\delta_{q;a_1,\dots,a_r} =
\delta_{q;ba_1,\dots,ba_r}$.

\item If the $a_j$ are either all squares modulo $q$ or all nonsquares
modulo $q$, then $\delta_{q;a_1,\dots,a_r} =
\delta_{q;a_r,\dots,a_1}$.

\item If $b$ is a reduced residue class modulo $q$ such that $c(q,a_j)\ne
c(q,ba_j)$ for each $1\le j\le r$, then $\delta_{q;a_1,\dots,a_r} =
\delta_{q;ba_r,\dots,ba_1}$.  In particular, this holds if $q$ is an
odd prime power or twice an odd prime power and $b$ is any nonsquare
modulo $q$.
\end{enumerate}
\end{theorem}

The pairwise equalities in Theorem 1 are special cases of part (d) of
Theorem 2, which generalizes the previously mentioned result of
Rubinstein and Sarnak that $\delta_{q;a,a'} = \delta_{q;a',a}$ if $a$
and $a'$ are either both squares or both nonsquares modulo $q$. Their
other symmetry result, that $\delta_{q;a_1,a_2,a_3}$ is invariant
under permutations of the $a_j$ when there exists $\rho \not\equiv
1\pmod q$ with $\rho^3 \equiv 1 \pmod q$ such that $a_2 \equiv
{a_1\rho} \pmod q$ and $a_3 \equiv {a_1\rho^2} \pmod q$, is also a
consequence of Theorem 2 (specifically parts (b) and (d), the former
applied with $b=\rho$ and $b=\rho^2$).

To complement Theorem 2, we can also establish several inequalities
concerning the densities $\delta$:

\begin{theorem}
Assume GRH and LI. Let $q\ge2$ be an integer,
let $N$ and $N'$ be distinct (invertible) nonsquares mod $q$, and let
$S$ and $S'$ be distinct (invertible) squares mod $q$. Then:
\begin{enumerate}
\item $\delta_{q;N,N',S} > \delta_{q;S,N',N}$;
\item $\delta_{q;N,S,S'} > \delta_{q;S',S,N}$;
\item $\delta_{q;N,S,N'} > \delta_{q;N',S,N}$ if and only if
$\delta_{q;N,S} > \delta_{q;N',S}$;
\item $\delta_{q;S,N,S'} > \delta_{q;S',N,S}$ if and only if
$\delta_{q;S,N} > \delta_{q;S',N}$.
\end{enumerate}
\end{theorem}

\noindent Parts (c) and (d) of Theorem 3 are further examples that the
predisposition towards some orderings of $\{\pi(x;q,a_1),\dots,\pi(x;q,a_r)\}$
over others cannot be explained solely in terms of the Chebyshev bias that encourages
nonsquares to run ahead of squares in the prime number race. (See also the
discussion of ``bias factors'' in Section 6.)

The most general result in this paper is an explicit formula for an
arbitrary density $\delta_\qaa$. Because of the amount of notation
involved, we have deferred the statement of this result (Theorem 4) to
Section 2.5. We have used this general formula to calculate the
densities given in Theorem 1, and also a number of the $\delta_\qaa$
in many interesting cases involving $q\le12$ and $r\le4$. In these
additional computations we have not undertaken to rigorously bound the
error terms; nevertheless we believe, from numerical considerations,
that the results given in Section 4 are accurate to the number
of decimal places indicated.

We shall assume the hypotheses GRH and LI throughout this paper. In
Section 2 we provide our main analysis leading to Theorem 4, the
general formula for $\delta_\qaa$.  The rigorous bounding of the error
terms incurred during the calculation of the densities in Theorem 1 is
carried out in Section 3.  Details of the computations and the
additional numerical results are collected together in Section 4.  The
proofs of Theorems 2 and 3 are given in Section 5, while in Section 6
we provide concluding remarks, noting some possible directions for
further work.

\bigskip
\section{Analytic Determination of the Densities $\delta_\qaa$.}
\setcounter{equation}{0}

The goal for this section of the paper is to derive Theorem 4 (see
Section 2.5), a general formula for the densities $\delta_\qaa$. We
begin by developing some notation and citing the relevant results of
Rubinstein and Sarnak in Section 2.1. In Section 2.2 we investigate
the function $\hat\rho_\qaa$ which will figure prominently in the
arguments that follow, while in Section 2.3 we establish some facts
about Cauchy principal values of multidimensional integrals; these sections
are technical rather than conceptual in nature,
and the reader may wish to examine these only briefly on the first reading.
Because the general formula
given in Theorem 4 and the arguments leading to it are somewhat involved, in
Section 2.4 we first detail the derivation of this formula for the special
cases $\delta_{8;a,b,c}$ and $\delta_{12;a,b,c}$ occurring in Theorem 1; the
derivation of the formula in the general case is then carried out in Section
2.5. We assume the hypotheses GRH and LI throughout.

\subsection{Notation and Background Results.}

We begin by establishing the notation necessary for discussing the
results of Rubinstein and Sarnak. For any coprime integers $q$ and $a$
and any real number $x\ge1$, define
\begin{equation}
E(x;q,a) = {\log x\over\sqrt x} \big( \phi(q)\pi(x;q,a)-\pi(x) \big),
\label{Exqadef}
\end{equation}
so that $E(x;q,a)$ is an error term for the number of primes congruent
to $a\bmod q$, normalized so as to vary roughly boundedly as $x$
varies. Since the inequalities $\pi(x;q,a_1)>\dots>\pi(x;q,a_r)$ hold if and
only if $E(x;q,a_1)>\dots>E(x;q,a_r)$, we wish to study how often the
vector
\begin{equation}
E_\qaa(x) = \big ( E(x;q,a_1),\dots,E(x;q,a_r) \big )
\label{Eqaadef}
\end{equation}
lies in the region $\{(x_1,\dots,x_r)\in\R^r\colon x_1>\dots>x_r\}$. Notice
that if $r=\phi(q)$ then the $a_j$ form a complete set of reduced residues
mod $q$, in which case we see from equation~(\ref{Exqadef}) that
\begin{equation}
E(x;q,a_1)+\dots+E(x;q,a_r)={-}{\log x\over\sqrt x}\phi(q)\omega(q)
\label{Esumequation}
\end{equation}
where $\omega(q)$ denotes the number of distinct prime factors of $q$.

Rubinstein and Sarnak showed, assuming GRH, that the function
$E_\qaa(x)$ has a limiting distribution $\mu_\qaa$, in the sense that
\begin{equation}
\lim_{X \rightarrow \infty} \frac{1}{\log X} \int_2^X
f(E_\qaa(x)) \, \frac{dx}{x} = \int\dots\int_{\R^r} f(x_1,\dots,x_r) \,
d\mu_\qaa
\label{limitingdistribution}
\end{equation}
for all bounded, continuous functions $f$ on $\R^r$. Under the further
assumption of LI, they showed that the distribution $\mu_\qaa$ is
absolutely continuous with respect to the ordinary Lebesgue measure on
$\R^r$. (The exception is the case $r=\phi(q)$, when
equation~(\ref{Esumequation}) implies that the distribution $\mu_\qaa$
is supported on the hyperplane $x_1+\dots+x_r=0$; in this case,
$\mu_\qaa$ is absolutely continuous with respect to Lebesgue measure
on this hyperplane.)  Consequently, the
equation~(\ref{limitingdistribution}) holds when $f$ is the
characteristic function of any reasonable subset of $\R^r$
(specifically, a measurable subset whose boundary has Lebesgue measure
zero in $\R^r$). In particular, it follows from the definition of
$\delta_\qaa$ that
\begin{equation}
\begin{split}
\delta_\qaa &= \delta \big ( \{ x \in \R\colon \pi(x,q,a_1) > \dots >
\pi(x,q,a_r) \} \big ) \\
&= \mu_\qaa \big ( \{ x \in \R^r\colon x_1 > \dots > x_r \} \big ) \\
&= \mathop{\int \dots \int}_{x_1 > \dots
> x_r} d \mu_\qaa \, .
\end{split}
\label{deltaasmu}
\end{equation}
Another consequence of the absolute continuity of $\mu_\qaa$ is that the
set of positive real numbers $x$ for which $\pi(x;q,a)=\pi(x;q,a')$ has
density zero when $a$ and $a'$ are distinct reduced residues;
indeed this is even true of the larger set $\{x\colon
|\pi(x;q,a)-\pi(x;q,a')| < \Phi(x)\}$ for any function $\Phi$ such
that
\begin{equation*}
\lim_{x\to\infty} {\Phi(x)\over\sqrt x/\log x}=0.
\end{equation*}

Next we develop the notation needed to write down Rubinstein and
Sarnak's seminal formula for the Fourier transform $\hat\mu_\qaa$ of
the distribution $\mu_\qaa$. In this paper we use the normalization
\begin{equation}
\hat f(\xi_1,\dots,\xi_n) = \int\dots\int e^{-i(\xi_1x_1 + \dots +
\xi_nx_n)} f(x_1,\dots,x_n) \, dx_1\dots dx_n
\label{normalization}
\end{equation}
for the Fourier transform of an integrable function $f$ on $\R^n$, so that
the Fourier inversion formula (assuming that $\hat f$ is itself integrable) is
\begin{equation*}
f(x_1,\dots,x_n) = (2\pi)^{-n}\int\dots\int e^{i(\xi_1x_1 + \dots +
\xi_nx_n)} \hat f(\xi_1,\dots,\xi_n) \, d\xi_1 \dots d\xi_n \, .
\end{equation*}
Likewise we write
\begin{equation*}
\hat \mu(\xi_1,\dots,\xi_n) = \int\dots\int e^{-i(\xi_1x_1 + \dots +
\xi_nx_n)} \, d\mu
\end{equation*}
for the Fourier transform of a finite measure $\mu$ on $\R^r$, so
that the Fourier inversion formula (assuming that $\hat \mu$ is integrable
with respect to Lebesgue measure) is
\begin{equation}
d\mu = (2\pi)^{-n} \bigg( \int\dots\int e^{i(\xi_1x_1 + \dots +
\xi_nx_n)} \hat\mu(\xi_1,\dots,\xi_n) \, d\xi_1 \dots d\xi_n \bigg) \,
dx_1 \dots dx_n.
\label{muinversion}
\end{equation}

To write down the specific Fourier transform $\hat\mu_\qaa$, we recall
the standard Bessel function of order zero,
\begin{equation}
J_0(z) = \sum_{m=0}^\infty \frac {(-1)^m (z/2)^{2m}} {(m!)^2}
= 1 - \frac{z^2}{4} + \frac {z^4} {64} - \cdots,
\label{Besseldef}
\end{equation}
and then set
\begin{equation}
F(z,\chi) = \prod \begin{Sb}
\gamma > 0 \\
L( \frac{1}{2} + i \gamma ,\chi) = 0
\end{Sb}
J_0(\alpha_\gamma z)
\label{Fzchidefn}
\end{equation}
in terms of the Dirichlet $L$-function $L(s,\chi)$ corresponding to the
Dirichlet character $\chi$, where we have defined
\begin{equation}
\alpha_\gamma = \frac {2} {\sqrt{ \frac{1}{4} + \gamma^2 }} \, .
\label{alphagammadefn}
\end{equation}
(Since we are assuming GRH, the product in equation~(\ref{Fzchidefn})
is indexed by all the nontrivial zeros of $L(s,\chi)$ in the upper
half-plane.) For later use in numerical approximations of $F(z,\chi)$
we also define the truncated version
\begin{equation}
F_T(z,\chi) = \bigg( \prod \begin{Sb}
0 < \gamma < T \\
L( \frac{1}{2} + i \gamma ,\chi) = 0
\end{Sb}
J_0(\alpha_\gamma z) \bigg) ( 1 + b_1 z^2 )
\label{FTzchidefn}
\end{equation}
for any positive real number $T$, where
\begin{equation}
b_1 = b_1(T,\chi) = - \sum_{\gamma \geq T} \frac {1} {\frac{1}{4} +
\gamma^2} \, .
\label{b1def}
\end{equation}
The polynomial factor in the definition~(\ref{FTzchidefn}) of $F_T$ is
motivated by the fact that, in view of the power series
expansion~(\ref{Besseldef}) of $J_0$, $b_1$ is the coefficient of
$z^2$ in the power series expansion of $\prod_{\gamma>T}
J_0(\alpha_\gamma z)$.

With this notation in place, we can now give the formula~\cite[equation
1.2]{RS} of Rubinstein and Sarnak for the Fourier transform $\hat\mu_\qaa$
of the distribution $\mu_\qaa$. They showed,
assuming GRH and LI, that
\begin{equation}
\hat\mu_\qaa (\xi_1, \dots , \xi_r) =
\exp \bigg( i \sum_{j=1}^r c(q,a_j)\xi_j \bigg) \prod
\begin{Sb} \chi \bmod q \\ \chi\ne\chi_0 \end{Sb}  F\bigg( \bigg|
\sum_{j=1}^r \chi(a_j)\xi_j \bigg|, \chi \bigg),
\label{muhatformula}
\end{equation}
where $c(q,a)$ was defined in equation~(\ref{cqadef}). This result
will be used extensively in the sequel.

Since $J_0(0)=1$ we clearly have $F(0,\chi)=F_T(0,\chi)=1$ for any
character~$\chi$. It is known (see for instance the arguments
in~\cite[Chapters 15--16]{davenport}) that for a fixed character $\chi$, the
number of zeros of $L(s,\chi)$ with imaginary part between 0 and $T$ has
order of magnitude $T\log T$. From this it can be shown that the
product~(\ref{Fzchidefn}) defining $F(z,\chi)$ converges uniformly on
bounded subsets of the complex plane, and hence $F$ is an entire function. 
For later use we will need bounds for the decay rate of $F(x,\chi)$ and
its derivatives $F^{(N)}(x,\chi)$ on the real axis;
this is the subject of the following lemma.

\begin{lemma}
Given a modulus $q\ge2$ and a nonnegative integer $N$, there
exist positive constants $\beta_1$ and $\beta_2$ such that
\begin{equation*}
|F^{(N)}(x,\chi)| \le \beta_1 e^{-\beta_2|x|}
\end{equation*}
for all real numbers $x$.
\label{Fboundlemma}
\end{lemma}

\noindent We caution the reader that in the next three sections, the
constants $\beta_1$ and $\beta_2$ will not necessarily have the same values
at different occurrences; each statement should be interpreted as holding for
some suitable positive values of $\beta_1$ and $\beta_2$.

\medskip\noindent {\bf Proof:} 
In this proof we will use the symbol $\gamma$, with or without subscript,
exclusively to denote a positive imaginary part of a nontrivial zero of
$L(s,\chi)$. We also use $\Gamma$ to denote an ordered $N$-tuple
$(\gamma_1,\dots,\gamma_N)$, and we let $m_\Gamma(\gamma)$ denote the number
(possibly zero) of coordinates of $\Gamma$ that equal $\gamma$. When
convenient we can also assume that $x>1$, since $F$ is an even, smooth
function. From the definition~(\ref{Fzchidefn}) of $F(z,\chi)$, an $N$-fold
application of the product rule gives us the expression
\begin{equation}
\begin{split}
F^{(N)}(x,\chi) &= \sum_{\Gamma=(\gamma_1,\dots,\gamma_N)} \alpha_{\gamma_1}
\dots \alpha_{\gamma_N} \prod_\gamma J_0^{(m_\Gamma(\gamma))}(\alpha_\gamma
x) \\
&= \sum_\Gamma \Phi(x,\Gamma) F(x,\chi,\Gamma)
\end{split}
\label{productrule}
\end{equation}
for the $N$th derivative of $F(x,\chi)$, where we have set
\begin{equation}
\Phi(x,\Gamma)= \prod_{\gamma\in\Gamma} \alpha_\gamma^{m_\Gamma(\gamma)}
J_0^{(m_\Gamma(\gamma))}(\alpha_\gamma x)
\label{PhixGammadef}
\end{equation}
and
\begin{equation*}
F(x,\chi,\Gamma) = \prod_{\gamma\notin\Gamma} J_0(\alpha_\gamma x).
\end{equation*}

We can show that $F(x,\chi)$ decays rapidly on the real axis by using
the standard bound~\cite[equation~(4.5)]{RS}
\begin{equation*}
| J_0(x) | \le \min \left \{ 1, \sqrt { \frac {2} {\pi |x|} } \, \right \}
\end{equation*}
for the Bessel function on the real axis. This bound implies that
\begin{equation}
\begin{split}
|F(x,\chi)| &\le \prod_\gamma \min \left \{ 1, \sqrt { \frac
{2} {\pi |\alpha_\gamma x|} } \, \right \} \\
&\leq \prod_{j=1}^J \sqrt { \frac {2} {\pi | \alpha_{\gamma_j}x |} } =
\ (\pi|x|)^{-J/2} \prod_{j=1}^J \left ( \frac{1}{4} + \gamma_j^2
\right )^{1/4}
\end{split}
\label{FboundedwithJ}
\end{equation}
for any positive integer $J$, where the $\gamma_j$ have been indexed
in increasing order. Choose $J=J(x)$ to be the number of zeros of $L(s,\chi)$
up to height $x/2$. For any $0<\gamma\le x/2$, it is easily verified that
the factor $(\pi|x|)^{-1/2} (1/4+\gamma^2)^{1/4}$ does not exceed $1/2$.
Therefore the upper bound~(\ref{FboundedwithJ}) implies that
$|F(x,\chi)|\le2^{-J}$. Since the order of magnitude of $J$ is $x\log x$,
this argument shows that as $x$ tends to infinity, $|F(x,\chi)|$ decreases
at least as fast as a function of the form $c^{x\log x}$ for some constant
$c$ depending on~$\chi$.

The same conclusion holds for $F(x,\chi,\Gamma)$, since
removing the indices $j$ in equation~(\ref{FboundedwithJ}) for which
$\gamma_j\in\Gamma$ changes $J$ by at most $N$ and thus does not affect the
order of magnitude of~$J$. Certainly then there exist positive constants
$\beta_1$ and $\beta_2$ (depending only on $N$ and $\chi$) such that
$|F(x,\chi,\Gamma)| \le \beta_1 e^{-\beta_2|x|}$ for all real numbers $x$.
Since this implies from equation~(\ref{productrule}) that
\begin{equation}
|F(x,\chi)| \le \beta_1 e^{-\beta_2|x|} \sum_\Gamma \Phi(x,\Gamma),
\label{justneedPhi}
\end{equation}
the lemma will be established (possibly with different values of $\beta_1$
and $\beta_2$) if we can show that this last sum is bounded by some
polynomial function of $|x|$.

To this end, we employ the crude bounds $|J_0'(t)| \le
\frac{t}{2}$ and $|J_0^{(n)}(t)| \le 1$ for the derivatives of the
Bessel function, which follow easily from the integral representation
\begin{equation*}
J_0(t) = {2\over\pi} \int_0^{\pi/2} \cos(t\sin\theta) \,d\theta.
\end{equation*}
Again supposing that $x>1$, the definition~(\ref{PhixGammadef}) of
$\Phi(x,\Gamma)$ leads to the bound
\begin{equation*}
|\Phi(x,\Gamma)| \le \bigg( \prod \begin{Sb}\gamma\in\Gamma \\
m_\Gamma(\gamma)=1\end{Sb} \alpha_\gamma^2 |x| \bigg)
\bigg( \prod \begin{Sb}\gamma\in\Gamma \\ m_\Gamma(\gamma)>1\end{Sb}
\alpha_\gamma^{m_\Gamma(\gamma)} \bigg).
\end{equation*}
It follows that
\begin{equation*}
\begin{split}
\bigg| \sum_\Gamma \Phi(x,\Gamma) \bigg| &\le |x|^N \sum_\Gamma
\alpha_\gamma^{\max\{m_\Gamma(\gamma),2\}} \\
&\le |x|^N N! \prod_\gamma (1 + \alpha_\gamma^2 + \alpha_\gamma^3 + \dots +
\alpha_\gamma^N).
\end{split}
\end{equation*}
Since the $j$th constant $\alpha_\gamma$ has order of magnitude $1/\gamma_j
\sim (\log j)/j$, this last product converges to some constant depending only
on $\chi$. Combining this bound with the inequality~(\ref{justneedPhi})
establishes the lemma.
\qed\bigskip


Of course it also follows from the first line of
equation~(\ref{FboundedwithJ}) that $|F(x,\chi)|$ is bounded above by 1
on the real axis.

In Sections 3.1 and 3.5 we will need to make use of the fact that $\mu_\qaa$
can also be thought of as the joint distribution of a certain set of
$r$ real-valued random variables, and it is convenient to exhibit
these random variables explicitly at this time. For given values of
$q$, $r$, and $a_1, \dots , a_r$, define the vector
\begin{equation*}
b_\qaa  =  {-} \big ( c(q,a_1), \dots ,  c(q,a_r) \big ) .
\end{equation*}
Next, for any character $\chi \bmod q$, define both the vector
\begin{equation*}
v_\qaa(\chi) = \big ( \chi(a_1), \dots , \chi(a_r) \big )
\end{equation*}
and the random variable
\begin{equation}
X(\chi) = \sum \begin{Sb} \gamma > 0 \\ L( \frac{1}{2} + i\gamma, \chi
) = 0 \end{Sb} \alpha_\gamma \sin ( 2\pi U_\gamma ),
\label{Xchi}
\end{equation}
where the $\alpha_\gamma$ are as in~(\ref{alphagammadefn}) and the
$U_\gamma$ are independent random variables uniformly distributed on
$[0,1]$. Note that by the hypothesis LI, the $\gamma$'s corresponding to
different $L$-functions are distinct, so that a given $U_\gamma$ only
appears in the definition of one of the $X(\chi)$; consequently the
random variables $ \{ X(\chi) \}$ are mutually independent. Then
Rubinstein and Sarnak showed that the distribution $\mu_\qaa$ is in fact
the same as the probability measure corresponding to the random vector
\begin{equation}
b_\qaa  +   \sum \begin{Sb}
\chi \bmod q \\
\chi \neq \chi_0
\end{Sb}
X(\chi) v_\qaa(\chi).
\label{bXchi}
\end{equation}

\subsection{The function $\hat\rho_\qaa$.}
In this section we introduce the function $\hat\rho_\qaa\colon
\R^{r-1}\to{\bf C}$, which we define by the formula
\begin{equation}
\hat\rho_\qaa(\eta_1,\dots,\eta_{r-1}) =
\hat\mu_\qaa(\eta_1,\eta_2-\eta_1,\dots,\eta_{r-1}-\eta_{r-2},-\eta_{r-1}),
\label{rhohatdef}
\end{equation}
so that
\begin{multline}
\hat\rho_\qaa(\eta_1,\dots,\eta_{r-1}) = \exp\bigg( \sum_{j=1}^{r-1}
\big ( c(q,a_j)-c(q,a_{j+1}) \big ) \eta_j \bigg) \\
\times \prod \begin{Sb}\chi\bmod q \\ \chi\ne\chi_0\end{Sb} F\bigg(
\bigg| \sum_{j=1}^{r-1} \big ( \chi(a_j)-\chi(a_{j+1}) \big ) \eta_j \bigg|, \chi
\bigg)
\label{rhohatformula}
\end{multline}
from the formula~(\ref{muhatformula}) for $\hat\mu_\qaa$. We will see
in Sections 2.4 and  2.5 that $\hat\rho_\qaa$ is the Fourier transform of a certain
measure $\rho_\qaa$ on $\R^{r-1}$ associated with $\mu_\qaa$.
We remark that in the special case
where the $a_j$ are all squares or all nonsquares, we have
$c(q,a_1)=\dots=c(q,a_r)$ and so the exponential term in the
formula~(\ref{rhohatformula}) is identically 1, so that $\hat\rho_\qaa$
is real-valued and symmetric with respect to reflection through the origin.

The function $\hat\rho_\qaa$ will feature significantly in the
remainder of this paper, and it will be important to establish some of
its smoothness and decay properties. To avoid frequent repetition of
the same properties, we shall say that a function $f$ on $\R^n$ is
{\it well-behaved\/} if it has continuous derivatives of all orders
and if there exist positive constants $\beta_1$ and $\beta_2$ such
that, for every subset $\{j_1,\dots,j_k\}$ of $\{1,\dots,n\}$, the
mixed partial derivative ${\partial^k f\over\partial x_{j_1} \dots
\partial x_{j_k}}$ satisfies the inequality
\begin{equation}
\left | {\partial^k f\over\partial x_{j_1} \dots \partial x_{j_k}}
(x_1,\dots,x_n) \right | \le \beta_1 e^{-\beta_2\|x\|},
\label{wellbineq}
\end{equation}
where $\|x\| = \|(x_1,\dots,x_n)\| = \sqrt{x_1^2+\dots+x_n^2}$ is the
Euclidean norm of $x$. This criterion must also be satisfied for the
empty subset of $\{1,\dots,n\}$, so that the actual values of $f$ must
also be bounded by the right-hand side of~(\ref{wellbineq}).
Certainly any well-behaved function is integrable as well.  We remark
that all of the functions shown to be well-behaved below in fact
satisfy an inequality analogous to~(\ref{wellbineq}) for partial
derivatives of all orders; however our proof of Lemma 2.4 below only
requires this assumption on the mixed linear partial derivatives.

It is easily seen that finite sums and products of well-behaved
functions are again well-behaved. If $f$ and $g$ are well-behaved
functions on $\R^m$ and $\R^n$, respectively, then $fg$ is a
well-behaved function on $\R^{m+n}$; conversely, the restriction of a
well-behaved function on $\R^n$ to any subspace defined by setting
certain variables equal to zero is a well-behaved function on that
subspace. Also, if $L\colon \R^m\to\R^n$ is an injective linear map
and $f$ is a well-behaved function on $\R^n$, then the composite
function $f\circ L$ is a well-behaved function on $\R^m$: the partial
derivatives of $f\circ L$ will just be linear combinations of the
partial derivatives of $f$, and the fact that $L$ is injective means
that $\|L(x)\|$ is bounded below by a constant multiple of $\|x\|$, so
that the estimate~(\ref{wellbineq}) for $f$ on $\R^n$ can be converted
to a similar estimate for $f\circ L$ on $\R^m$.

The following two lemmas establish the important fact that the functions
$\hat\rho_\qaa$ are well-behaved.

\begin{lemma}
For every subset $\{j_1,\dots,j_k\}$ of $\{1,\dots,r\}$,
\begin{equation}
{\partial^k \over\partial x_{j_1} \dots \partial x_{j_k}} F\bigg( \bigg|
\sum_{j=1}^r \chi(a_j)x_j \bigg|, \chi \bigg) \le \beta_1 \exp \bigg(
{-}\beta_2 \bigg| \sum_{j=1}^r \chi(a_j)x_j \bigg| \bigg)
\label{harderFbound}
\end{equation}
for some positive constants $\beta_1$ and $\beta_2$.
\label{harderFboundlemma}
\end{lemma}

\noindent {\bf Proof:}
Define
\begin{equation*}
G(x_1,\dots,x_r;\chi) = F\bigg( \bigg| \sum_{j=1}^r \chi(a_j)x_j \bigg|,
\chi \bigg).
\end{equation*}
The argument of $F$ on the right-hand side of this definition involves
a modulus and hence implicitly a square root, which could potentially cause
discontinuities in the derivatives of $G$ when this argument equals zero;
however, the Bessel function $J_0$ is even, whence the function $F(x,\chi)$
involves only even powers of $x$ in its power series expansion about the
origin. Consequently, $G$ has continuous derivatives of all orders.
Note also that it suffices to establish the upper bound~(\ref{harderFbound})
when $|\sum_{j=1}^r \chi(a_j)x_j|>1$, since the bound
on the complementary set follows immediately
from continuity (with some value of $\beta_1$).

If we write $\tilde F(x,\chi) = F(\sqrt{|x|},\chi)$, then it is easy to check
by induction that the $n$th derivative of $\tilde F$ equals
\begin{equation*}
\tilde F^{(n)}(x,\chi) = \sum_{k=1}^n \alpha_{n,k} F^{(k)}
\big( \sqrt{|x|},\chi \big) |x|^{-n+j/2}
\end{equation*}
for some constants $\alpha_{n,k}$. In particular, when $|x|>1$ we see from
Lemma~\ref{Fboundlemma} that
\begin{equation}
|\tilde F^{(n)}(x,\chi)| \le \beta_1 e^{{-}\beta_2 \sqrt{|x|}}
\label{Ftildebound}
\end{equation}
for some positive constants $\beta_1$ and $\beta_2$.

In this notation we have
\begin{equation*}
G(x_1,\dots,x_r;\chi) = \tilde F\bigg( \bigg( \Re \sum_{j=1}^r \chi(a_j)x_j
\bigg)^2 + \bigg( \Im \sum_{j=1}^r \chi(a_j)x_j \bigg)^2, \chi \bigg).
\end{equation*}
Suppressing the details, we note that the mixed partial derivative
${\partial^k G\over\partial x_{j_1} \dots \partial x_{j_k}}$ can be
computed using the product rule as a combination of three types of
expressions: derivatives of $\tilde F$ evaluated at $|\sum_{j=1}^r
\chi(a_j)x_j|^2$, linear factors of the form $2\Re( \bar\chi(a_k)
\sum_{j=1}^r \chi(a_j)x_j )$, and constants of the form
$2\Re(\bar\chi(a_k)\chi(a_{k'}))$.  From equation~(\ref{Ftildebound}),
the expressions of the first type can be bounded above by $\beta_1
\exp ( {-}\beta_2 | \sum_{j=1}^r \chi(a_j)x_j | )$, while the
expressions of the other types grow only as fast as a polynomial in $|
\sum_{j=1}^r \chi(a_j)x_j |$. This establishes the lemma for suitable
positive values of $\beta_1$ and $\beta_2$.\qed\bigskip

\begin{lemma}
The function $\hat\rho_\qaa$ is well-behaved for any integers $q,r\ge2$ and
any distinct reduced residues $\{a_1,\dots,a_r\}$.
\label{rhohatWBlemma}
\end{lemma}

\noindent{\bf Proof:} From the formula~(\ref{muhatformula}), the function
$\hat\mu_\qaa$ certainly has continuous derivatives of all orders (see the
proof of Lemma~\ref{harderFboundlemma}), and thus the same is true of
$\hat\rho_\qaa$. We begin by examining the behavior of the mixed partial
derivatives of the function $\hat\mu_\qaa$. Let $S=\{j_1,\dots,j_k\}$ be a subset
of indices from the set $\{1,\dots,r\}$, and let ${\partial^k\over\partial
x_S}$ denote the result of taking the partial $x_j$-derivatives for every $j$
in $S$. The product rule applied to
the formula~(\ref{muhatformula}) for $\hat\mu_\qaa$ yields
\begin{equation}
{\partial^k \over \partial x_S} \hat\mu_\qaa(\xi) = \sum_{S_0,\{S_\chi\}}
\bigg\{ {\partial^k \over \partial x_{S_0}} \exp \bigg( i \sum_{j=1}^r
c(q,a_j)\xi_j \bigg) \prod \begin{Sb} \chi \bmod q \\ \chi\ne\chi_0
\end{Sb} {\partial^k \over \partial x_{S_\chi}} F\bigg( \bigg|
\sum_{j=1}^r \chi(a_j)\xi_j \bigg|, \chi \bigg) \bigg\},
\label{productruletomuhat}
\end{equation}
where the outer summation is taken over the finitely many partitions
of the index set $S$ into $S_0 \cup (\bigcup_{\chi\ne\chi_0}
S_\chi)$. Each mixed partial derivative of the exponential term is
bounded, while from Lemma~\ref{harderFboundlemma} each mixed partial
derivative of $F(|\sum_{j=1}^r \chi(a_j)\xi_j|,\chi)$ is exponentially
decaying as a function of its argument. We conclude from
equation~(\ref{productruletomuhat}) that there exist positive
constants $\beta_1$ and $\beta_2$ such that
\begin{equation}
\big| {\partial^k \over \partial x_S} \hat\mu_\qaa(\xi) \big| \le
\beta_1 \prod \begin{Sb}\chi\bmod q \\ \chi\ne\chi_0\end{Sb} \exp \bigg(
{-\beta_2} \bigg| \sum_{j=1}^r \chi(a_r)\xi_j \bigg| \bigg) =
\beta_1 e^{{-\beta_2} Q(\xi)^{1/2}},
\label{Qintroduced}
\end{equation}
where we have defined
\begin{equation*}
Q(\xi) = Q_\qaa(\xi) = \bigg( \sum \begin{Sb}\chi\bmod q \\
\chi\ne\chi_0\end{Sb} \bigg| \sum_{j=1}^r \chi(a_j)\xi_j \bigg|
\bigg)^2.
\end{equation*}
We thus seek a lower bound on $Q(\xi)$.

We may certainly write
\begin{equation*}
Q(\xi) \ge \sum \begin{Sb}\chi\bmod q \\ \chi\ne\chi_0\end{Sb} \bigg|
\sum_{j=1}^r \chi(a_j)\xi_j \bigg|^2 = \sum_{\chi\bmod q} \bigg|
\sum_{j=1}^r \chi(a_j)\xi_j \bigg|^2 - \bigg( \sum_{j=1}^r \xi_j
\bigg)^2.
\end{equation*}
Now
\begin{equation*}
\sum_{\chi\bmod q} \bigg|
\sum_{j=1}^r \chi(a_j)\xi_j \bigg|^2 = \sum_{i=1}^r \sum_{j=1}^r \xi_i
\xi_j \sum_{\chi\bmod q} \chi(a_i)\overline{\chi(a_j)} = \phi(q)
\sum_{j=1}^r \xi_j^2
\end{equation*}
by the orthogonality of the characters $\chi$. Therefore
\begin{equation}
Q(\xi) \ge \phi(q) \sum_{j=1}^r \xi_j^2 - \bigg( \sum_{j=1}^r \xi_j
\bigg)^2.
\label{Qxilowerbound}
\end{equation}

We assume for now that $r$ is strictly less than $\phi(q)$, commenting
at the end of the proof on the slight differences in the case
$r=\phi(q)$. The quadratic form on the right-hand side of the
inequality~(\ref{Qxilowerbound}) turns out to be positive definite
when $r<\phi(q)$, and so we can write
\begin{equation}
Q(\xi) \ge \phi(q)\lambda_r \|\xi\|^2,  \label{strangeQ}
\end{equation}
where $\lambda_r$ is the smallest eigenvalue of that quadratic form. From the
inequalities~(\ref{Qintroduced}) and~(\ref{strangeQ}), it follows that
\begin{equation*}
\big| {\partial^k \over \partial x_S} \hat\mu_\qaa(\xi) \big| \le
\beta_1 e^{{-\beta_2}\|\xi\|}
\end{equation*}
for some different positive constants $\beta_1$ and $\beta_2$. Since
the index set $S\subset\{1,\dots,r\}$ was arbitrary, this shows that
the function $\hat\mu_\qaa$ is well-behaved.

Furthermore, from its definition~(\ref{rhohatdef}) the function
$\hat\rho_\qaa$ is simply the composition of $\hat\mu_\qaa$ with the
injective linear transformation $(\eta_1,\dots,\eta_{r-1}) \mapsto
(\eta_1,\eta_2-\eta_1,\dots,\eta_{r-1}-\eta_{r-2},-\eta_{r-1})$ from
$\R^{r-1}$ to $\R^r$. As mentioned before, this implies that $\hat\rho_\qaa$
is itself a well-behaved function.

When $r=\phi(q)$, the function $\hat\mu_\qaa$ is invariant under translation
in the direction of the vector $(1,\dots,1)$, and so it is not well-behaved
even though it has the required decay properties on the hyperplane orthogonal
to $(1,\dots,1)$ (one can check that the quadratic form on the right-hand
side of the inequality~(\ref{Qxilowerbound}) is positive semi-definite when
$r=\phi(q)$, with its zero set being the multiples of the $(1,\dots,1)$ vector). However,
the image of the linear transformation $(\eta_1,\dots,\eta_{r-1}) \mapsto
(\eta_1,\eta_2-\eta_1,\dots,\eta_{r-1}-\eta_{r-2},-\eta_{r-1})$ lies within
this hyperplane, so we can still deduce that $\hat\rho_\qaa$ is well-behaved
even when $r=\phi(q)$. This establishes the lemma.
\qed\bigskip

Of course we also have the trivial bound $|\hat\rho_\qaa|\le1$.
Lemma~\ref{rhohatWBlemma} implies in particular that $\hat\rho_\qaa$ is
integrable, and consequently the Fourier inversion formula~(\ref{muinversion})
is valid for $\rho_\qaa$, becoming
\begin{equation}
d\rho_\qaa = (2\pi)^{-r} \bigg( \int\dots\int e^{i(\xi_1x_1 + \dots +
\xi_rx_r)} \hat\rho_\qaa(\xi_1,\dots,\xi_r) \, d\xi_1 \dots d\xi_r
\bigg) \, dx_1 \dots dx_r.
\label{rhoinversion}
\end{equation}

\subsection{Multidimensional Cauchy principal values.}
In one dimension, the Cauchy principal value
\begin{equation*}
\PV \int_{-\infty}^\infty {f(x)\over x} \, dx = \lim_{\epsilon\to0}
\int_{|x|>\epsilon} {f(x)\over x} \, dx
\end{equation*}
is a familiar object. For our purposes it will be necessary to make
use of the multidimensional analogue
\begin{equation}
\PV \int \dots \int {f(x_1,\dots,x_n)\over x_1\dots x_n} \, dx_1 \dots
dx_n = \lim_{\epsilon\to0}
\mathop{\int\dots\int}_{\min\{|x_1|,\dots,|x_n|\}>\epsilon}
{f(x_1,\dots,x_n)\over x_1\dots x_n} \, dx_1 \dots dx_n;
\label{multiPV}
\end{equation}
in particular, we would like to know that this limit exists. The purpose of
this section is to establish the existence of these multidimensional Cauchy
principal values for well-behaved functions, a class which by
Lemma~\ref{rhohatWBlemma} includes the functions $\hat\rho_\qaa$ discussed in
the previous section. We remark that while the lemmas in
this section could certainly be obtained under somewhat weaker hypotheses,
they suffice for our purposes as stated.

\begin{lemma}
Let $f$ be a well-behaved function on $\R^n$ that vanishes
whenever any of the first $k$ coordinates $x_1$, \dots, $x_k$ equals zero.
Then the function $f(x_1,\dots,x_n)/x_1\dots x_k$ extends across the
coordinate hyperplanes to a continuous integrable function satisfying
the upper bound
\begin{equation}
\big| {f(x_1,\dots,x_n)\over x_1\dots x_k} \big| \le
\beta_1e^{{-\beta_2}\|x\|}
\label{dividingOKbound}
\end{equation}
for some positive constants $\beta_1$ and $\beta_2$.
\label{dividingOKlemma}
\end{lemma}

\noindent Although this lemma holds in one dimension without any
assumptions on the derivatives of $f$, already in $\R^2$ one can construct an
exponentially decaying, smooth (even real-analytic) function $f(x,y)$ that
satisfies $f(0,y)=0$ for all $y$ but for which $f(x,y)/x$ is not integrable.

\medskip\noindent {\bf Proof:} The fact that $f(x)/x_1\dots x_k$ extends
across the coordinate hyperplanes to a continuous function follows from the
fact that $f$ has continuous derivatives of all orders; therefore only
the upper bound~(\ref{dividingOKbound}) remains to be proved, since
integrability is a consequence of this bound. Furthermore, by continuity
it suffices to establish this upper bound when none of the
variables equals zero. Also, if all of the $|x_j|$ are bounded by 1
then the function $f(x)/x_1\dots x_k$ is uniformly bounded; therefore we may
assume (after inflating the constant $\beta_1$ if necessary) that there exists
an $x_j$ with $|x_j|>1$.

Permuting the first $k$ variables if necessary, we can choose an integer $1\le
m\le k$ such that $0<|x_1|,\dots,|x_m|\le1$ and $|x_{m+1}|,\dots,|x_k|>1$.
Since $f$ vanishes when $x_1$ equals zero, there exists a number $t_1$ with
$|t_1|\le|x_1|$ such that
\begin{equation*}
f(x_1,\dots,x_n) = f(x_1,\dots,x_n) - f(0,x_2,\dots,x_n) = x_1{\partial
f\over\partial x_1}(t_1,x_2,\dots,x_n)
\end{equation*}
by the mean value theorem in the variable $x_1$. Similarly, $f$ vanishes
whenever $x_2$ equals zero, so in particular ${\partial f\over\partial x_1}$
equals zero when $x_2=0$. Therefore, there exists a number $t_2$ with
$|t_2|\le|x_2|$ such that
\begin{equation*}
\begin{split}
{\partial f\over\partial x_1}(t_1,x_2,\dots,x_n) &= {\partial f\over\partial
x_1}(t_1,x_2,\dots,x_n) - {\partial f\over\partial x_1}(t_1,0,x_3,\dots,x_n)
\\
&= x_2{\partial^2 f\over\partial x_1\partial x_2}(t_1,t_2,x_3,\dots,x_n)
\end{split}
\end{equation*}
by the mean value theorem in the variable $x_2$. Continuing in this way, we
find numbers $t_i$ with $|t_i|\le|x_i|$ for each $1\le i\le m$ such that
\begin{equation*}
f(x_1,\dots,x_n) = x_1\dots x_m{\partial^m f\over\partial x_1\dots\partial
x_m}(t_1,\dots,t_m,x_{m+1},\dots,x_n).
\end{equation*}
It follows immediately that
\begin{equation}
\big| {f(x)\over x_1\dots x_k} \big| \le \big| {\partial^m f\over\partial
x_1\dots\partial x_m}(t_1,\dots,t_m,x_{m+1},\dots,x_n) \big|
\label{meanvalues}
\end{equation}
since $|x_{m+1}|,\dots,|x_k|>1$.

Since $f$ is well-behaved, there exist positive constants $\beta_1$
and $\beta_2$ such that
\begin{equation}
\big| {\partial^m f\over\partial x_1\dots\partial x_m} \big(
t_1,\dots,t_m,x_{m+1},\dots,x_n) \big| \le \beta_1\exp\big (
{-\beta_2} \sqrt{t_1^2 + \dots + t_m^2 + x_{m+1}^2 + \dots + x_n^2}
\big).
\label{useWBwithts}
\end{equation}
But notice that
\begin{equation*}
t_1^2 + \dots + t_m^2 + x_{m+1}^2 + \dots + x_n^2 \ge \sum \begin{Sb}1\le
j\le n \\ |x_j|>1\end{Sb} x_j^2 \ge {\#\{1\le j\le n\colon |x_j|>1\}\over n}
\sum_{j=1}^n x_j^2.
\end{equation*}
Since we are working under the assumption that at least one of the $|x_j|$
exceeds 1, we can use this fact in the inequality~(\ref{useWBwithts}) to see
that
\begin{equation*}
\big| {\partial^m f\over\partial x_1\dots\partial
x_m}(t_1,\dots,t_m,x_{m+1},\dots,x_n) \big| \le \beta_1 e^{{-\beta_2}
\|x\|/\sqrt n}.
\end{equation*}
Combining this bound with the inequality~(\ref{meanvalues}), this establishes
the lemma (upon replacing $\beta_2/\sqrt n$ by $\beta_2$).
\qed\bigskip

For the proof of the next lemma, as well as for the formulation of the general
formula for $\delta_\qaa$ (Theorem 4 in Section 2.5), we require the following
notation: for a function $f$ on $\R^n$ and a subset $B$ of $\{1,\dots,n\}$,
define
\begin{equation}
f(B) = f(B)( \{ x_j\colon j \in B \} ) = f (\theta_1, \dots, \theta_n)
\label{subsetnotation}
\end{equation}
where $\theta_j = x_j$ if $j \in B$, and $\theta_j = 0$ otherwise.
For example, if $n=6$ and $B = \{ 2,4,5 \}$ then $f(B)$ is a function
of the three variables $x_2$, $x_4$, and $x_5$, namely $f(B) = f(0,
x_2, 0, x_4, x_5,0)$; in general $f(B)$ will be a function on the
appropriate $|B|$-dimensional subspace of $\R^n$, where $|B|$ denotes the
cardinality of $B$. In the case $B=\emptyset$ we simply have
$f(B)=f(0,\dots,0)$.

\begin{lemma}
If $f$ is a well-behaved function on $\R^n$, then
\begin{equation*}
\PV \int \dots \int {f(x_1,\dots,x_n)\over x_1\dots x_n} \, dx_1 \dots dx_n
\end{equation*}
is well-defined; i.e., the limit in equation~(\ref{multiPV}) exists.
\label{PVexists}
\end{lemma}

\noindent{\bf Proof.}
Let $g_1(x)$ be an even, well-behaved function on $\R^1$ with $g_1(0)=1$ (for
instance, we might have in mind $g_1(x)=e^{-x^2}$), and let
$g(x_1,\dots,x_n)=g_1(x_1)\dots g_1(x_n)$. Define an operator $G$ on
well-behaved functions $f$ by
\begin{equation}
G(f) = G(f)(x_1,\dots,x_n) = \sum_{B \subset \{1,\dots,n\}}
(-1)^{n-|B|} f(B) g(\bar B)
\label{Gfdef}
\end{equation}
in the notation of equation~(\ref{subsetnotation}), where $\bar B$ denotes the
complement $\{1,\dots,n\}\setminus B$ of $B$. Since $f$ and $g$ are
well-behaved functions, the same is true of $G(f)$.

Consider the term in (\ref{Gfdef}) corresponding to some particular proper subset $B$
of $\{1,\dots,n\}$. If we choose $\ell\notin B$, then the term
$f(B)g(\bar B)$ can be written as $g(x_\ell)$ times a function
independent of $x_\ell$. Thus $f(B)g(\bar B)$ is an even function of
$x_\ell$, and hence integrates to zero against any odd function of
$x_\ell$. In particular,
\begin{equation*}
\mathop{\int\dots\int}_{\min\{|x_1|,\dots,|x_n|\}>\epsilon}
{f(B)g(\bar B) \over x_1\dots x_r} \, dx_1 \dots dx_n = 0
\end{equation*}
for any positive $\epsilon$ and any proper subset $B$ of $\{1,\dots,n\}$.
Since the term in the sum~(\ref{Gfdef}) corresponding to $B=\{1,\dots,n\}$ is
simply the function $f$ itself, we see that
\begin{equation}
\mathop{\int\dots\int}_{\min\{|x_1|,\dots,|x_n|\}>\epsilon}
{f(x_1,\dots,x_n)\over x_1\dots x_n} \, dx_1 \dots dx_n =
\mathop{\int\dots\int}_{\min\{|x_1|,\dots,|x_n|\}>\epsilon} {G(f)\over
x_1\dots x_n} \, dx_1 \dots dx_n
\label{Gtakesover}
\end{equation}
for any $\epsilon>0$.

On the other hand, we claim that $G(f)$ evaluates to zero when any
of the variables $x_\ell$ equals zero. To see this, let $B$ be a
subset of $\{1,\dots,n\}$ not containing $\ell$. When $x_\ell=0$ we
see that the term $(-1)^{n-|B|} f(B)g(\bar B)$ corresponding to $B$ in
the sum~(\ref{Gfdef}) reduces to $(-1)^{n-|B|} f(B)g(\bar
B\setminus\{\ell\})$. On the other hand, the term
\begin{equation*}
(-1)^{n-|B\cup\{\ell\}|} f(B\cup\{\ell\})g(\overline{B\cup\{\ell\}})
= (-1)^{n-1-|B|} f(B\cup\{\ell\})g(\bar B\setminus\{\ell\})
\end{equation*}
corresponding to $B\cup\{\ell\}$ reduces to $(-1)^{n-1-|B|} f(B)g(\bar
B\setminus\{\ell\})$ when $x_\ell = 0$. It follows that when $x_\ell = 0$,
the terms in~(\ref{Gfdef}) will cancel pairwise in the natural
pairing between the subsets of $\{1,\dots,n\}$ not containing $\ell$ and
those containing $\ell$.

Because of this, Lemma~\ref{dividingOKlemma} tells us that the function
$G(f)/x_1\dots x_n$ is integrable, whence the dominated convergence theorem
implies
\begin{equation}
\lim_{\epsilon\to0}
\mathop{\int\dots\int}_{\min\{|x_1|,\dots,|x_n|\}>\epsilon} {G(f)\over
x_1\dots x_n} \, dx_1 \dots dx_n = \int\dots\int {G(f)\over
x_1\dots x_n} \, dx_1 \dots dx_n.
\label{qwert}
\end{equation}
This together with equation~(\ref{Gtakesover}) shows that the principal
value~(\ref{multiPV}) exists---in fact it equals the integral on the
right-hand side of equation~\ref{qwert}.\qed\bigskip

\begin{lemma}
If $f$ is a well-behaved function on $\R^n$, then for any $1\le k\le n$,
\begin{multline}
\lim_{c\to0+} c^{n-k} \int\dots\int_{\R^n}
{f(x_1,\dots,x_n) x_1\dots x_k \over (c^2+x_1^2) \dots (c^2+x_n^2)}
\, dx_1\dots dx_n \\
= \pi^{n-k}\, \PV \int\dots\int_{\R^k} {f(x_1,\dots,x_k,0,\dots,0) \over
x_1\dots x_k} \, dx_1\dots dx_k.
\label{limitlemmaeqn}
\end{multline}
\label{limitlemma}
\end{lemma}

\medskip
\noindent{\bf Proof.}
We proceed along lines similar to the proof of Lemma~\ref{PVexists}.
Analogously to the definition~(\ref{Gfdef}) of the operator $G(f)$, define
the operator
\begin{equation*}
G_k(f) = G_k(f)(x_1,\dots,x_n) = \sum_{B \subset \{1,\dots,k\}}
(-1)^{k-|B|} f(B\cup\{k+1,\dots,n\}) g(\bar B),
\end{equation*}
so that $G_k(f)$ is itself a well-behaved function. The arguments leading to
the validity of equation~(\ref{Gtakesover}) in the proof of
Lemma~\ref{PVexists} show that the function $G_k(f)-f$ integrates to 0
against any function that is odd in each of the variables $x_1$, \dots, $x_k$
separately. In particular,
\begin{multline*}
c^{n-k} \int\dots\int {f(x_1,\dots,x_n) x_1\dots x_k \over
(c^2+x_1^2) \dots (c^2+x_n^2)} \, dx_1\dots dx_n \\
= c^{n-k} \int\dots\int {G_k(f)(x_1,\dots,x_n) x_1\dots x_k \over
(c^2+x_1^2) \dots (c^2+x_n^2)} \, dx_1\dots dx_n.
\end{multline*}
Making the change of variables $x_j = c\xi_j$ for $k < j \leq n$ and
rearranging terms, we see that
\begin{multline}
c^{n-k} \int\dots\int {G_k(f)(x_1,\dots,x_n) x_1\dots x_k \over
(c^2+x_1^2) \dots (c^2+x_n^2)} \, dx_1\dots dx_n \\
= \int\dots\int  {\tilde G_k(f)(x_1,\dots,x_k,cx_{k+1},\dots,cx_n) x_1^2\dots
x_k^2 \over (c^2+x_1^2) \dots (c^2+x_k^2) (1+x_{k+1}^2)\dots(1+x_n^2)}
\, dx_1\dots dx_n,
\label{needsupfn}
\end{multline}
where we have defined
\begin{equation*}
\tilde G_k(f)(x) = {G_k(f)(x)\over x_1\dots x_k}.
\end{equation*}
As in the proof of Lemma~\ref{PVexists}, we can check that $G_k(f)$ evaluates
to zero whenever any of the first $k$ variables equals zero, and thus by
Lemma~\ref{dividingOKlemma} the function $\tilde G_k(f)$ is
continuous and integrable and satisfies an upper bound of the form
\begin{equation}
|\tilde G_k(f)(x)| \le \beta_1e^{{-\beta_2}\|x\|}  \label{twostar}
\end{equation}
for some positive constants $\beta_1$ and $\beta_2$.

Now define
\begin{equation*}
S_c(x_1,\dots,x_k) = \begin{cases}
\beta_1e^{-\beta_2c\|x\|}, &\hbox{if }\|x\| > 1/\sqrt c, \\
\max\limits_{|t_{k+1}|, \dots, |t_n| \le \sqrt c} {\displaystyle \left
| {\tilde G_k(f)(x_1,\dots,x_k,t_{k+1},\dots,t_n) \over
(1+x_{k+1}^2)\dots(1+x_n^2)} \right | }, &\hbox{if }\|x\| \le 1/\sqrt c.
\end{cases}
\end{equation*}
One can check that the integrand on the right-hand side of
equation~(\ref{needsupfn}) is bounded in absolute value by
$S_c(x_1,\dots,x_k)$ when $0<c<1$. Moreover, the continuity of $\tilde
G_k(f)$ implies that $S_c$ is bounded on the set $\{x\in\R^k\colon
\|x\|\le 1/\sqrt c\}$, and therefore $S_c$ is integrable. Furthermore,
both $S_c$ and the integrand on the right-hand side of
equation~(\ref{needsupfn}) tend pointwise to the function
\begin{equation*}
{\tilde G_k(f)(x_1,\dots,x_k,0,\dots,0) \over (1+x_{k+1}^2)\dots(1+x_n^2)}
\end{equation*}
as $c$ tends to zero, and this function is also integrable by the
exponential decay~\ref{twostar} of $\tilde G_k(f)$. Therefore, taking
limits on both sides of equation~(\ref{needsupfn}) and using the
generalized dominated convergence theorem, we conclude that
\begin{equation*}
\begin{split}
\lim_{c\to0+} c^{n-k} \int\dots\int & {G_k(f)(x_1,\dots,x_n) x_1\dots x_k
\over (c^2+x_1^2) \dots (c^2+x_n^2)} \, dx_1\dots dx_n \\
&= \int\dots\int {\tilde G_k(f)(x_1,\dots,x_k,0,\dots,0) \over
(1+x_{k+1}^2)\dots(1+x_n^2)} \, dx_1\dots dx_n \\
&= \pi^{n-k} \int\dots\int {G_k(f)(x_1,\dots,x_k,0,\dots,0) \over x_1\dots
x_k} \, dx_1\dots dx_k.
\end{split}
\end{equation*}
But just as in the proof of Lemma~\ref{PVexists}, this last integral equals
the principal value of the integral of $f(x_1,\dots,x_k,0,\dots,0)/x_1\dots
x_k$, which establishes the lemma.\qed\bigskip

\noindent Of course, the lemma would also hold if both occurrences of the
product $x_1\dots x_k$ in equation~(\ref{limitlemmaeqn})
were replaced by any product $x_{j_1}\dots x_{j_k}$ of
$k$ distinct variables (and the variables of integration on the right-hand
side adjusted accordingly).

\subsection{Analysis for the special case.}
In this section we derive analytic expressions for the density
$\delta_\qnn$ (the logarithmic density of the set $\{ x \in
\R\colon \pi(x;8,3) > \pi(x;8,5) > \pi(x;8,7) \}$) and for the other
densities in Theorem 1. These formulas are special cases of Theorem~4,
which will be established in the next section; however, we present a
complete analysis in these special cases to
illustrate and motivate the techniques in the general case.

We begin by noting the special case
\begin{equation*}
\delta_\qnn = \mathop{\int \int \int}_{x>y>z} d \mu_\qnn (x,y,z)
\end{equation*}
of equation~(\ref{deltaasmu}).
Making the change of variables $u=x-y$, $v=y-z$, and $w=z$ gives
\begin{equation*}
\delta_\qnn = \int\limits_{u>0} \int\limits_{v>0} \int\limits_{w\in\R} d
\nu_\qnn(u,v,w),
\end{equation*}
where the measure $\nu_\qaa$ is defined, in obvious notation, by
\begin{equation}
\nu_\qnn (u,v,w) = \mu_\qnn ( u+v+w, v+w, w ),
\label{mu2nu}
\end{equation}
or equivalently $\mu_\qnn (x,y,z) = \nu_\qnn (x-y, y-z, z)$.
Integrating out the $w$ variable, we obtain
\begin{equation}
\delta_\qnn = \int\limits_{u>0} \int\limits_{v>0} d \rho_\qnn (u,v),
\label{deltarhohat}
\end{equation}
where we have defined, again in obvious notation,
\begin{equation}
\rho_\qnn (u,v) = \int\limits_{w\in\R} d \nu_\qnn (u,v,w) \, .
\label{nu2rho}
\end{equation}
It is easily checked that the Fourier transform of $\rho_\qnn$ is
related to that of $\mu_\qnn$ via $\hat\rho_\qnn(\xi,\eta) =
\hat\mu_\qnn(\xi,\eta-\xi,-\eta)$, which is a particular case of
equation~(\ref{rhohatdef}).

We can appeal to the formula~(\ref{rhohatformula}) for $\hat\rho_\qaa$
to write $\hat\rho_\qnn(\xi,\eta)$ explicitly. Recall that a
discriminant is an integer congruent to 0 or 1 mod 4, and a
fundamental discriminant $D$ is an integer that cannot be written in
the form $D=dn^2$ for some discriminant $d$ and integer $n\ge2$. For
any fundamental discriminant $D$, let $\chi_D$ denote the character
$\chi_D(n)=({D\over n})$ using Kronecker's extension of the Legendre
symbol (see Davenport~\cite[Chapter 5]{davenport}). Then the three
nonprincipal characters mod 8 are simply $\chi_{-8}$, $\chi_{-4}$, and
$\chi_8$. In this setting, equation~(\ref{rhohatformula}) becomes
\begin{equation}
\hat\rho_\qnn(\xi,\eta) = F(|2\xi|,\chi_{-8}) F(|2\eta-2\xi|,\chi_{-4})
F(|-2\eta|,\chi_8),
\label{rhohatqnnformula}
\end{equation}
showing that $\hat\rho_\qnn$ is real-valued and symmetric with respect
to reflection through the origin. The same argument gives formulas for
$\delta_{q;a,b,c}$ for any permutation $\{a,b,c\}$ of $\{3,5,7\}$,
where the arguments of the $F(\cdot,\chi)$ functions in
equation~(\ref{rhohatqnnformula}) simply are permuted accordingly. Since each
$F(z,\chi)$ is an even function, we can omit the absolute value signs
in these arguments. Similar remarks hold for the modulus 12, where
the three nonprincipal characters are $\chi_{-4}$, $\chi_{-3}$, and
$\chi_{12}$.

Using the monotone convergence theorem and the Fourier inversion
formula~(\ref{rhoinversion}), equation~(\ref{deltarhohat}) becomes
\begin{equation*}
\begin{split}
\delta_\qnn &= \lim_{c \rightarrow 0+} \int\limits_{u>0} \int\limits_{v>0}
e^{-c(u+v)} \, d \rho_\qnn (u,v) \\
&= \lim_{c \rightarrow 0+} \int\limits_{u>0} \int\limits_{v>0} e^{-c(u+v)}
\left \{ \frac{1}{4 \pi^2} \int \int e^{i(u\xi+v\eta)} \hat \rho_\qnn
(\xi,\eta) \, d\xi d\eta \right \} \, du dv.
\end{split}
\end{equation*}
We next use Fubini's theorem to write
\begin{equation}
\begin{split}
\delta_\qnn &= \frac{1}{4 \pi^2} \lim_{c \rightarrow 0+} \mathop{\int
\int} \hat \rho_\qnn (\xi,\eta) \bigg\{ \, \int\limits_{u>0}
\int\limits_{v>0} e^{u(-c+i\xi) + v(-c+i\eta)} \, du dv \bigg\} \, d\xi
d\eta \\
&= \frac{1}{4 \pi^2} \lim_{c \rightarrow 0+} \mathop{\int \int} \frac
{\hat \rho_\qnn (\xi,\eta)} {(c-i\xi) (c-i\eta)} \, d\xi d\eta \\
&= \frac{1}{4 \pi^2} \lim_{c \rightarrow 0+} \mathop{\int \int} \frac
{\hat \rho_\qnn (\xi,\eta) (c^2 + i c(\xi+\eta) - \xi\eta)}
{(c^2+\xi^2) (c^2+\eta^2)} \, d\xi d\eta \\
&= \frac{1}{4 \pi^2} (G_\qnn +i H_\qnn - I_\qnn),
\label{GHI}
\end{split}
\end{equation}
where we have defined
\begin{equation}
G_\qnn = \lim_{c \rightarrow 0+} c^2 \mathop{\int \int} \frac {\hat
\rho_\qnn (\xi,\eta)} {(c^2+\xi^2) (c^2+\eta^2)} \, d\xi d\eta \, ,
\label{G}
\end{equation}
\begin{equation}
H_\qnn = \lim_{c \rightarrow 0+} c \mathop{\int \int} \frac {\hat
\rho_\qnn (\xi,\eta) (\xi+\eta) } {(c^2+\xi^2) (c^2+\eta^2)} \, d\xi
d\eta \, ,
\label{H}
\end{equation}
and
\begin{equation}
I_\qnn = \lim_{c \rightarrow 0+} \mathop{\int \int} \frac {\hat
\rho_\qnn (\xi,\eta) \xi \eta } {(c^2+\xi^2) (c^2+\eta^2)} \, d\xi
d\eta \, .
\label{I}
\end{equation}

In equation~(\ref{G}) we make the change of variables $\alpha = \xi/c$ and
$\beta = \eta/c$ to obtain
\begin{equation*}
G_\qnn = \lim_{c \rightarrow 0+} \mathop{\int \int} \frac {\hat
\rho_\qnn (c\alpha,c\beta)} {(1+\alpha^2) (1+\beta^2)} \, d\xi d\eta =
\mathop{\int \int} \frac {\hat \rho_\qnn (0,0)} {(1+\alpha^2)
(1+\beta^2)} \, d\alpha d\beta = \pi^2
\end{equation*}
where we have again used the dominated convergence theorem together
with the trivial bound $|\hat \rho_\qnn(\xi,\eta)| \le
\hat\rho_\qnn(0,0)=1$. Next, we note that $H_\qnn$ equals zero since
the integrand in equation~(\ref{H}) is odd under reflection through
the origin.  Finally, we observe that equation~(\ref{I}) may be
written as
\begin{equation*}
I_\qnn = \lim_{c \rightarrow 0+} \mathop{\int \int} \frac { \left (
\hat \rho_\qnn (\xi,\eta) - \hat \rho_\qnn (\xi,0) \hat \rho_\qnn
(0,\eta) \right ) \xi \eta } {(c^2+\xi^2) (c^2+\eta^2)} \, d\xi d\eta,
\end{equation*}
since the term introduced is odd in either variable separately
and so integrates to zero. This is the same as
\begin{equation*}
I_\qnn = \lim_{c \rightarrow 0+} \mathop{\int \int}
\frac { \hat \rho_\qnn (\xi,\eta) - \hat \rho_\qnn (\xi,0) \hat
\rho_\qnn (0,\eta) }  { \xi \eta}
\  \frac {\xi^2 \eta^2}  {(c^2+\xi^2) (c^2+\eta^2)} \  d\xi d\eta.
\end{equation*}

Note that the expression $\hat \rho_\qnn (\xi,\eta) - \hat \rho_\qnn (\xi,0)
\hat \rho_\qnn (0,\eta)$ is well-behaved by Lemma~\ref{rhohatWBlemma}, and
since $\hat\rho_\qnn(0,0)=1$, it evaluates to zero when either $\xi$ or
$\eta$ equals zero. Therefore, the first fraction in the integrand can be
extended across the $\xi$ and $\eta$ axes to a continuous integrable function
by Lemma~\ref{dividingOKlemma}. We may thus use the dominated convergence
theorem to see that
\begin{equation}
I_\qnn = \mathop{\int \int} \frac { \hat \rho_\qnn (\xi,\eta) - \hat
\rho_\qnn (\xi,0) \hat \rho_\qnn (0,\eta) } { \xi \eta}   \  d\xi d\eta.
\label{I8357}
\end{equation}
Note that this integral may be written as the multivariate Cauchy
principal value
\begin{equation}
I_\qnn = \PV \mathop{\int \int}
\frac { \hat \rho_\qnn (\xi,\eta) }  { \xi \eta}   \, d\xi d\eta
\label{I8357PV}
\end{equation}
as discussed in Section 2.3, since $\hat \rho_\qnn (\xi,0)$ and $\hat
\rho_\qnn (0,\eta)$ are even functions and hence the term omitted in passing
from (\ref{I8357}) to (\ref{I8357PV}) is odd in either variable. (Of course, we could
have arrived at (\ref{I8357PV}) directly from the definition of $I_\qnn$ by
invoking Lemma~\ref{limitlemma}; however, not only is this derivation
more concrete, in keeping with the spirit of this section, but we will also need the
formula~(\ref{I8357}) during our error analysis in Section 3.)

It follows that the right-hand side of equation~(\ref{GHI}) can be evaluated
to give
\begin{equation}
\delta_\qnn = \frac{1}{4} - \frac{1}{4 \pi^2} I_\qnn = \frac{1}{4} -
\frac{1}{4 \pi^2} \, \PV \mathop{\int \int} \frac { \hat \rho_\qnn
(\xi,\eta) } { \xi \eta} \, d\xi d\eta \, ,
\label{the8357formula}
\end{equation}
where $\hat\rho_\qnn$ is given explicitly in
equation~(\ref{rhohatqnnformula}). The identical argument, of course,
applies for evaluating $\delta_{8;a,b,c}$ for any permutation $\{a,b,c\}$ of
$\{3,5,7\}$ to yield
\begin{equation}
\delta_{8;a,b,c} \ =
\ \frac{1}{4} \, - \, \frac{1}{4 \pi^2} \, \PV \mathop{\int \int}
\frac { \hat \rho_{8;a,b,c} (\xi,\eta) }  { \xi \eta}   \, d\xi d\eta,
\label{obvanal}
\end{equation}
where $\rho_{8;a,b,c}$ is defined via obvious analogy to (\ref{mu2nu})
and (\ref{nu2rho}), and similarly
\begin{equation}
\delta_{12;a,b,c} \ =
\ \frac{1}{4} \, - \, \frac{1}{4 \pi^2} \, \PV \mathop{\int \int}
\frac { \hat \rho_{12;a,b,c} (\xi,\eta) }  { \xi \eta}   \, d\xi d\eta,
\label{obvanal12}
\end{equation}
for any permutation $\{a,b,c\}$ of $\{5,7,11\}$.

We remark that the numerator of the integrand in (\ref{I8357}) may be
viewed as a ``measure of dependence'' in the Fourier domain for the
bivariate distribution of a random vector $(X,Y)$ in $\R^2$ having
density $\rho_\qnn$.  In fact, the integrand in equation~(\ref{I8357})
is the Fourier transform of the natural dependence measure based on
factorizability of the bivariate cumulative distribution function
corresponding to $\rho_\qnn$. This interpretation is important in
Section 3.1, where a random vector $(X,Y)$ of this type is analyzed to
yield bounds for the tail of the measure $\rho_\qnn$.

\subsection{Analysis for the general case.}
We are now at the point where we have the notation and tools needed for
the statement and proof of a general formula for the densities $\delta_\qaa$.

\begin{theorem}
Assume GRH and LI. Let $q,r\ge2$ be integers, and let $a_1$, \dots, $a_r$ be
distinct reduced residue classes mod $q$. Then
\begin{equation}
\delta_\qaa \ = \ 2^{-(r-1)} \bigg( 1 + \sum \begin{Sb} B \subset \{ 1, 
\dots , r-1\} \\ B\ne\emptyset \end{Sb} \left ( \frac {i}{\pi} \right
)^{|B|} \, \PV \int\dots\int\hat\rho_\qaa (B) \prod_{j \in B} \frac { d
\eta_j } {\eta_j} \bigg),
\label{generalformula}
\end{equation}
where $\hat\rho_\qaa(B)$ uses the notation of
equation~(\ref{subsetnotation}) applied to the function
\begin{multline*}
\hat\rho_\qaa(\eta_1,\dots,\eta_{r-1}) = \exp\bigg( \sum_{j=1}^{r-1}
(c(q,a_j)-c(q,a_{j+1}))\eta_j \bigg) \\
\times \prod \begin{Sb}\chi\bmod q \\ \chi\ne\chi_0\end{Sb} F\bigg(
\bigg| \sum_{j=1}^{r-1} (\chi(a_j)-\chi(a_{j-1}))\eta_j \bigg|, \chi
\bigg).
\end{multline*}
\end{theorem}

\medskip\noindent{\bf Proof:}
We follow the strategy used for the special cases in Section 2.4. Our
starting point is equation~(\ref{deltaasmu}),
\begin{equation*}
\delta_\qaa = \mathop{\int \dots \int}_{x_1 > x_2 > \dots > x_r}
d \mu_\qaa (x_1, \dots , x_r) \, .
\end{equation*}
We make the change of variables $u_1=x_1-x_2,$ \dots, $u_{r-1}=x_{r-1}-x_r$,
$u_r=x_r$ to obtain
\begin{equation*}
\delta_\qaa = \mathop{\int \dots \int} \begin{Sb}
u_1 > 0,\, \dots ,\, u_{r-1} > 0 \\ u_r \in \R\end{Sb}
d \nu_\qaa (u_1, \dots , u_r) \, ,
\end{equation*}
where the measure $\nu_\qaa$ is defined, in obvious notation, by
\begin{equation*}
\nu_\qaa (u_1, \dots , u_r) =
\mu_\qaa ( u_1 + \dots + u_r, u_2 + \dots + u_r, \dots , u_r ) \, ,
\end{equation*}
or equivalently $\mu_\qaa (x_1, \dots , x_r) = \nu_\qaa ( x_1-x_2,
\dots , x_{r-1}-x_r, x_r )$. Integrating out the $u_r$ variable leads to
\begin{equation}
\delta_\qaa = \int\limits_{u_1>0} \dots \int\limits_{u_{r-1}>0} d \rho_\qaa
(u_1, \dots , u_{r-1}) \, ,
\label{integratedout}
\end{equation}
where we have defined (again in obvious notation)
\begin{equation*}
\rho_\qaa (u_1, \dots , u_{r-1}) =
\int\limits_{v \in \R} d \nu_\qaa (u_1, \dots , u_{r-1}, v) \, .
\end{equation*}
It is easily checked that the Fourier transform of $\rho_\qaa$ is
related to that of $\mu_\qaa$ by the identity~(\ref{rhohatdef}).

At this point, our goal is to evaluate the integral on the right-hand
side of equation~(\ref{integratedout}) in terms of the Fourier
transform $\hat\rho_\qaa$ of $\rho_\qaa$. The correct final formula
could be obtained by writing this as the integral of $d\rho_\qaa$
against the characteristic function of the region of integration, and
using Parseval's identity in the context of the theory of generalized
functions; the following analysis derives this final formula
rigorously.

Using the monotone convergence theorem followed by the Fourier inversion
formula~(\ref{rhoinversion}), equation~(\ref{integratedout}) becomes
\begin{equation*}
\begin{split}
\delta_\qaa &= \lim_{c\to0+} \int\limits_{u_1>0} \dots \int\limits_{u_{r-1}>0}
e^{-c(u_1 + \dots + u_{r-1})} \, d \rho_\qaa (u_1, \dots , u_{r-1}) \\
&= \lim_{c\to0+} \int\limits_{u_1>0} \dots \int\limits_{u_{r-1}>0} e^{-c(u_1 +
\dots + u_{r-1})} \bigg\{ (2\pi)^{-(r-1)} \\
&\qquad\times \int \dots \int e^{i(u_1 \xi_1 + \dots + u_{r-1} \xi_{r-1})}
\hat \rho_\qaa (\xi_1, \dots , \xi_{r-1}) \, d\xi_1 \dots d\xi_{r-1} \bigg\}
\, du_1 \dots du_{r-1}.
\end{split}
\end{equation*}
Then by Fubini's theorem this becomes
\begin{equation*}
\begin{split}
\delta_\qaa &= (2\pi)^{-(r-1)} \lim_{c\to0+} \int \dots \int \hat \rho_\qaa
(\xi_1, \dots , \xi_{r-1}) \\
&\qquad\qquad\times \bigg\{ \int\limits_{u_1>0} \dots \int\limits_{u_{r-1}>0}
e^{u_1(-c+i\xi_1) + \dots + u_{r-1}(-c+i\xi_{r-1})} \, du_1
\dots du_{r-1} \bigg\} \, d\xi_1 \dots d\xi_{r-1} \\
&= (2\pi)^{-(r-1)} \lim_{c\to0+} \int \dots \int \frac {\hat \rho_\qaa (\xi_1,
\dots , \xi_{r-1})} {(c-i\xi_1)\dots(c-i\xi_{r-1})} \, d\xi_1 \dots
d\xi_{r-1} \\
&= (2\pi)^{-(r-1)} \lim_{c\to0+} \int \dots \int \frac {\hat \rho_\qaa
(\xi_1, \dots , \xi_{r-1}) (c+i\xi_1)\dots(c+i\xi_{r-1}) }
{(c^2+\xi_1^2)\dots(c^2+\xi_{r-1}^2)} \, d\xi_1 \dots d\xi_{r-1},
\end{split}
\end{equation*}
and expanding the product $(c+i\xi_1)\dots(c+i\xi_{r-1})$ leads to
\begin{equation}
\delta_\qaa = (2\pi)^{-(r-1)} \sum_{ B \subset \{ 1,2,...,{r-1} \} } i^{|B|}
\, I_\qaa (B),
\label{subsetexpand}
\end{equation}
where we have defined
\begin{equation*}
I_\qaa (B) = \lim_{c\to0+} c^{{r-1}-|B|} \int \dots \int \frac {\hat \rho_\qaa
(\xi_1, \dots , \xi_{r-1}) \left ( \prod_{j \in B} \xi_j \right ) }
{(c^2+\xi_1^2)\dots(c^2+\xi_{r-1}^2)} \, d\xi_1 \dots d\xi_{r-1} \, .
\end{equation*}
Appealing to Lemma~\ref{limitlemma} with $n=r-1$ and $k=|B|$, we see that
\begin{equation*}
I_\qaa(B) = \pi^{{r-1}-|B|} \, \PV \int \dots \int \hat \rho_\qaa (B) \prod_{j
\in B} {d\eta_j\over\eta_j},
\end{equation*}
which includes the special case $I_\qaa(\emptyset)=\pi^{r-1}$. Using this fact
in equation~(\ref{subsetexpand}) establishes the theorem.\qed\bigskip

We remark that the measure $\rho_\qaa$ is actually the limiting distribution
of the vector
\begin{equation*}
{\log x\over\sqrt x} \big( \pi(x;q,a_1)-\pi(x;q,a_2), \dots,
\pi(x;q,a_{r-1})-\pi(x;q,a_r) \big)
\end{equation*}
in $\R^{r-1}$, so its usefulness to the investigation of those $x$ with
$\pi(x;q,a_1)>\dots>\pi(x;q,a_r)$ is not surprising.

To conclude this section, we consider two special cases of Theorem 4.
In the case $r=2$ (in other words, when we are comparing simply a
pair $a_1,a_2$ of residues modulo $q$) the formula
(\ref{generalformula}) reduces to
\begin{equation}
\begin{split}
\delta_{q;a_1,a_2} &= \frac12 \bigg( 1 + \frac i\pi \, \PV \int
{\hat\rho_{q;a_1,a_2}(\eta)d\eta\over\eta} \bigg)  \\
&= \frac12 - {1\over2\pi} \int
{\sin( \{ c(q,a_1)-c(q,a_2) \} \eta)\over\eta} \prod \begin{Sb}\chi\bmod q \\
\chi\ne\chi_0\end{Sb} F_\chi(|\chi(a_1)-\chi(a_2)|\eta) \,d\eta \, ,
\end{split}
\label{genRS}
\end{equation}
the corresponding cosine term in the last integral being
omitted by virtue  of symmetry. When $c(q,a_1)=c(q,a_2)$, the
integrand is identically zero and hence $\delta_{q;a_1,a_2}=1/2$, as
was proved by Rubinstein and Sarnak. In fact, the
formula~(\ref{genRS}) is analogous to their formula~\cite[equation
4.1]{RS}.

In the case $r=3$, Theorem 4 becomes
\begin{multline}
\delta_{q;a_1,a_2,a_3} = \frac14 + \frac i{4\pi} \, \PV \int
\left ( \hat\rho_{q;a_1,a_2,a_3}(\eta,0) + \hat\rho_{q;a_1,a_2,a_3}(0,\eta)
\right ) \, {d\eta\over\eta} \\
- {1\over4\pi^2} \, \PV \int\int \hat\rho_{q;a_1,a_2,a_3}(\eta_1,\eta_2)
\, {d\eta_1 d\eta_2\over\eta_1 \eta_2}.
\end{multline}
If the $a_j$ are all squares or all nonsquares, then the one-dimensional
integral again vanishes due to symmetry, yielding a generalization of the
formulas~(\ref{obvanal}) and~(\ref{obvanal12}) of Section~2.4.

\bigskip
\section{Rigorous Error Bounds.}
\setcounter{equation}{0}

In this section, we describe how the densities in Theorem 1 were calculated
and provide a rigorous analysis bounding the error between the
calculated and true values.

Suppose that we wish to evaluate $\delta_\qnn$. According to
equation~(\ref{the8357formula}), we need only to evaluate
\begin{equation}
\PV \, \mathop{\int \int} \frac { \hat \rho_\qnn (\xi,\eta) } {
\xi \eta} \, d\xi d\eta = \PV \, \mathop{\int \int} \frac {
F(2\xi,\chi_{-8}) F(2\eta-2\xi,\chi_{-4}) F(-2\eta,\chi_8) } { \xi
\eta} \, d\xi d\eta,
\label{wanted}
\end{equation}
where we have used the formula~(\ref{rhohatqnnformula}) for
$\hat\rho_\qnn$. We shall approximate this integral by sampling the
integrand on the (symmetrically offset) grid of points
\begin{equation*}
\left \{
\left ( \frac{m\epsilon}{2} , \frac {n\epsilon}{2} \right ) : \,
\big| {m\epsilon\over2} \big|, \big| {n\epsilon\over2} \big| \leq C ,
\ m,n \ {\rm odd}
\right \}
\end{equation*}
for some appropriately small $\epsilon > 0$ and some appropriately large
$C>0$. In fact the quantity we actually compute is
$4 S_\qnn(\epsilon,C,T)$, where we define
\begin{equation}
S_\qnn(\epsilon,C,T) = \mathop{\sum \sum} \begin{Sb}|m|, |n| \leq
2C/\epsilon \\ m,n{\rm\ odd}\end{Sb}
\frac {F_T(m\epsilon,\chi_{-8}) F_T((n-m)\epsilon,\chi_{-4})
F_T(-n\epsilon,\chi_8)}  {mn};
\label{SCT}
\end{equation}
here $F_T(z,\chi)$ is the approximation to $F(z,\chi)$ defined in
equation~(\ref{FTzchidefn}), and as before $\chi_D$ is the character
given by the Kronecker symbol $\chi_D(n)=({D\over n})$.

The quantity $S_\qnn(\epsilon,C,T)$ is a discrete, truncated
approximation to the integral~(\ref{wanted}) involving an approximated
summand as well.  The overall error incurred in evaluating
(\ref{wanted}) by means of (\ref{SCT}) thus consists of three
components: error due to discretizing the integral, error due to
truncating the resulting infinite sum, and error due to approximating
the summand. In Sections~3.1--3.3, respectively, we obtain rigorous
bounds for each of these sources of error, and in Section~3.4 we
combine these bounds to establish Theorem~1. Section~3.5 provides some
technical bounds that are required for our arguments in
Section~3.1. While in the sections to follow, all of the specific
expressions we write down (such as $S_\qnn(\epsilon,C,T)$) are those
that arise in the calculation of the single density $\delta_\qnn$, the
given constants and error bounds were chosen so as to apply also to
the analogous quantities arising during the calculation of any of the
densities listed in Theorem~1.

\subsection{Error Due To Discretization.}

The first step is to discretize the calculation of $I_\qnn$ by converting
the integral defining $I_\qnn$ into a sum; we may bound the error incurred
by doing so using the Poisson summation formula, as we now explain. Let
$f(\xi,\eta)$ be a continuous, integrable function on $\R^2$ such that both
$f$ and $\hat f$ decay rapidly enough near infinity (for instance,
exponential decay certainly suffices). Then $f$ satisfies the Poisson
summation formula
\begin{equation*}
\epsilon_1 \epsilon_2 \sum_{k=-\infty}^\infty \sum_{\ell=-\infty}^\infty
f(k\epsilon_1 + \alpha, \ell \epsilon_2 + \beta )  =
\sum_{\kappa=-\infty}^\infty \sum_{\lambda=-\infty}^\infty
\hat f \left ( \frac{2\pi\kappa}{\epsilon_1} ,
\frac{2\pi\lambda}{\epsilon_2}
\right ) e^{2\pi i (\kappa \alpha / \epsilon_1 + \lambda \beta /
\epsilon_2)}
\end{equation*}
(see for instance Stein and Weiss~\cite[Corollary 2.6 of Chapter
VII]{SW}, although we are using a Fourier
transform~(\ref{normalization}) with a different choice of
constants). In this formula, set $\epsilon_1=\epsilon_2=\epsilon$ and
$\alpha = \beta \ = \epsilon / 2$, and make the change of variables $m
= 2k+1$ and $n = 2 \ell + 1$ on the left-hand side, to obtain
\begin{equation}
\epsilon^2 \, \mathop{\sum \sum} \begin{Sb}  m, n \in \Z \\  m, n {\rm\ 
odd} \end{Sb} f \left ( \frac {m\epsilon} {2} , \frac {n\epsilon} {2}
\right ) = \hat f(0,0) \ + \mathop{\sum \sum} \begin{Sb} \kappa, \lambda
\in\Z \\ (\kappa, \lambda)\neq (0,0) \end{Sb} \hat f \left ( \frac
{2\pi\kappa} {\epsilon} , \frac {2\pi\lambda} {\epsilon} \right )
(-1)^{\kappa + \lambda}.
\label{specialPoisson}
\end{equation}

Now let
\begin{equation*}
f(\xi,\eta) = \frac { \hat \rho_\qnn (\xi, \eta) - \hat \rho_\qnn
(\xi, 0) \hat \rho_\qnn (0 , \eta) } {\xi \eta},
\end{equation*}
which can be extended continuously over the coordinate axes as was
noted in Section 2.3. This function $f$ is integrable and has
exponential decay near infinity by Lemmas~\ref{rhohatWBlemma}
and~\ref{dividingOKlemma}, and its Fourier transform can be seen to
equal
\begin{equation}
\hat f (u,v) = -4\pi^2(\overline P (u,v) - \overline P_1(u) \overline
P_2(v)),
\label{fintermsofP}
\end{equation}
where
\begin{equation*}
\overline P (u,v) = \int_u^\infty \int_v^\infty d \rho_\qnn
\end{equation*}
is the upper cumulative distribution function of $\rho_\qnn$ and
$\overline P_1(u) = \overline P (u,-\infty)$ and $\overline P_2(v) =
\overline P (-\infty,v)$ are the corresponding ``upper marginals''. (Note
that $\hat f (u,v)$ is a dependence measure of the type mentioned at the
end of Section 2.4.) At the end of this section we will show that the
function $\hat f$ decays exponentially as well, so that we are
justified in applying the form~(\ref{specialPoisson}) of the
Poisson summation formula to~$f$.

Now observe from equation~(\ref{I8357}) that
\begin{equation*}
I_\qnn = \int \int f(\xi,\eta)\, d\xi d\eta = \hat f (0,0).
\end{equation*}
Therefore applying equation~(\ref{specialPoisson}) to the function
$f$, we have
\begin{equation}
I_\qnn = \epsilon^2 \mathop{\sum \sum} \begin{Sb} m, n \in \Z \\ m,
n {\rm\ odd} \end{Sb} \frac { \hat \rho_\qnn ( m\epsilon/2,
n\epsilon/2 ) - \hat \rho_\qnn ( m\epsilon/2, 0 ) \hat \rho_\qnn ( 0,
n\epsilon/2 ) } {(m\epsilon/2)(n\epsilon/2)} + {\rm Error}_1,
\label{I8Error}
\end{equation}
where Error$_1$, the error due to discretization, is given by
\begin{equation}
{\rm Error}_1  =  4\pi^2
\mathop{\sum \sum} \begin{Sb} \kappa, \lambda
\in\Z \\ (\kappa, \lambda)\neq (0,0) \end{Sb} \big( \, \overline P \big(
\frac{2\pi\kappa}{\epsilon},\frac{2\pi\lambda}{\epsilon} \big) -
\overline P_1 \big( \frac{2\pi\kappa}{\epsilon} \big)
\overline P_2 \big( \frac{2\pi\lambda}{\epsilon} \big)  \big)
(-1)^{\kappa + \lambda}.
\label{Error1}
\end{equation}
Defining
\begin{multline*}
Q(u,v) = \big( \overline P (2\pi u,2\pi v) - \overline P_1 (2\pi u)
\overline P_2 (2\pi v) \big) + \big( \overline P (-2\pi u,2\pi v) -
\overline P_1 (-2\pi u) \overline P_2 (2\pi v) \big) \\
+ \big( \overline P (2\pi u,-2\pi v) - \overline P_1 (2\pi u) \overline P_2
(-2\pi v) \big) + \big( \overline P (-2\pi u,-2\pi v) - \overline P_1 (-2\pi
u) \overline P_2 (-2\pi v) \big),
\end{multline*}
and grouping the terms on the right-hand side of equation~(\ref{Error1})
analogously, we obtain
\begin{equation*}
{\rm Error}_1 = 4\pi^2 \bigg( \mathop{\sum \sum}_{\kappa, \lambda \in
\Z^+} (-1)^{\kappa + \lambda} Q \big( \frac {\kappa} {\epsilon} ,
\frac {\lambda} {\epsilon} \big) + \frac{1}{2} \sum_{\kappa \in \Z^+}
(-1)^\kappa Q \big( \frac {\kappa} {\epsilon} , 0 \big) + \frac{1}{2}
\sum_{\lambda \in \Z^+} (-1)^\lambda Q \big( 0 , \frac {\lambda}
{\epsilon} \big) \bigg),
\end{equation*}
so that
\begin{equation}
| {\rm Error}_1 | \le 4\pi^2 \mathop{\sum \sum} \begin{Sb} \kappa,
\lambda \ge0 \\ (\kappa, \lambda)\neq (0,0) \end{Sb} \big| Q \big(
\frac {\kappa} {\epsilon} , \frac {\lambda} {\epsilon} \big) \big|.
\label{AbsError1}
\end{equation}

Now let $(X,Y)$ denote a pair of real-valued random variables whose joint
distribution is given by $\rho_\qnn$ (these random variables are given
explicitly in equation~(\ref{XandYdef}) below, though their explicit form
is not needed here).  Then $\overline P (u,v) = {\rm Pr} ( X>u,\, Y>v )$ and
hence $\overline P_1 (u) = {\rm Pr} ( X>u )$ and $\overline P_2 (v) = {\rm
Pr} ( Y>v )$.  With this interpretation, and using the fact that
$\rho_\qnn$ is symmetric about the origin, the identity
\begin{equation}
Q(u,v) = {\rm Pr} ( X>2\pi u, Y>2\pi v ) - {\rm Pr} ( X>2\pi u, Y<-2\pi v )
\label{Qidentity}
\end{equation}
is easily verified. Clearly
\begin{equation*}
0 \le {\rm Pr} ( X>u,\, Y>v ) \leq \min \{ {\rm Pr} ( X>u ) , {\rm Pr}
( Y>v ) \}.
\end{equation*}
Moreover, since $\rho_\qnn$ is symmetric about the origin,
each component $X$ and $Y$ is a symmetric random variable, so that
\begin{equation*}
0 \le {\rm Pr} ( X>u,\, Y<-v ) \leq \min \{ {\rm Pr} ( X>u ) , {\rm
Pr} ( Y<-v ) \} = \min \{ {\rm Pr} ( X>u ) , {\rm Pr} ( Y>v ) \}.
\end{equation*}
It therefore follows from the identity~(\ref{Qidentity}) that
\begin{equation}
| Q(u,v) | \le \min \{ {\rm Pr} ( X>2\pi u ) , {\rm Pr} ( Y>2\pi v ) \}.
\label{absQbound}
\end{equation}

In Section 3.5 we shall establish the bounds
\begin{equation}
\begin{split}
{\rm Pr} ( X \geq u ) &\le \exp ( {-0.04} (u-3)^2 ) \\
{\rm Pr} ( Y \geq u ) &\le \exp ( {-0.04} (u-3)^2 )
\end{split}
\label{thebounds}
\end{equation}
for any $u\ge3$. Hence by the inequality~(\ref{absQbound}),
\begin{equation*}
| Q(u,v) | \le \exp \big( {-0.04} (2\pi \max \{ u,v \} -3)^2 \big)
\end{equation*}
if either $u$ or $v$ exceeds 1, so that by equation~(\ref{AbsError1}),
\begin{equation}
\begin{split}
| {\rm Error}_1 | &\leq 4 \pi^2 \bigg ( \sum_{\kappa=0}^\infty
\sum_{\lambda = \max \{ \kappa, 1 \}}^\infty \exp \big( {-0.04} \big( \frac
{2\pi \lambda} {\epsilon} -3 \big)^2 \big) \\
&\qquad + \ \sum_{\lambda=0}^\infty \sum_{\kappa = \max \{ \lambda, 1
\}}^\infty \exp \big( {-0.04} \big( \frac
{2\pi \kappa} {\epsilon} -3 \big)^2 \big) \bigg) \\
&\leq \ 8 \pi^2 \sum_{\kappa=0}^\infty \sum_{\lambda = \max \{ \kappa, 1
\}}^\infty \exp \big( {-0.04} \big( \frac
{2\pi \lambda} {\epsilon} -3 \big)^2 \big)
\end{split}
\label{twice}
\end{equation}
if $\epsilon<1$, say. Now for any positive integer $\lambda_0$,
\begin{equation*}
\sum_{\lambda = \lambda_0}^\infty \exp \big( {-0.04} \big( \frac {2\pi
\lambda} {\epsilon} -3 \big)^2 \big) \leq 2 \exp \big( {-0.04} \big(
\frac {2\pi \lambda_0} {\epsilon} -3 \big)^2 \big),
\end{equation*}
since each term of the sum is at most half of the preceding term.
Applying this inequality twice in succession to the
bound~(\ref{twice}) gives
\begin{equation*}
\begin{split}
| {\rm Error}_1 | &\leq 8 \pi^2 \bigg( 2 \exp \big( {-0.04} \big(
\frac{2\pi}\epsilon -3 \big)^2 \big) + 2 \sum_{\kappa = 1}^\infty \exp \big(
{-0.04} \big( \frac{2\pi\kappa}\epsilon -3 \big)^2 \bigg) \\
&\leq 48 \pi^2 \exp \big( {-0.04} \big(
\frac{2\pi}\epsilon -3 \big)^2 \big).
\end{split}
\end{equation*}
We therefore conclude that
\begin{equation}
| {\rm Error}_1 | < 5\times 10^{-12}
\label{Error1bound}
\end{equation}
for any choice of $\epsilon<1/5$, which is more than adequate for our
purposes.

To conclude this section, we return to the matter of showing that the
function $\hat f$ given in equation~(\ref{fintermsofP}) decays
exponentially. In terms of the random variables $X$ and $Y$, the
formula~(\ref{fintermsofP}) becomes
\begin{equation*}
\hat f(u,v) = -4\pi^2 ( {\rm Pr} (X>u,\, Y>v) - {\rm Pr}(X>u) {\rm Pr}(Y>v) ).
\end{equation*}
By an argument similar to the one used for the function $Q$, we see that
\begin{equation}
\begin{split}
| \hat f(u,v) | &\le 4\pi^2 \max\{ {\rm Pr} (X>u,Y>v), {\rm Pr}(X>u) {\rm
Pr}(Y>v) \} \\
&\le 4\pi^2 \min\{ {\rm Pr}(X>u), {\rm Pr}(Y>v) \}.
\end{split}
\label{zazaza}
\end{equation}
On the other hand, by elementary considerations we have
\begin{equation*}
\begin{split}
\hat f(-u,-v) &= 4\pi^2 ( {\rm Pr} (X>-u,\, Y>-v) - {\rm Pr}(X>-u) {\rm
Pr}(Y>-v) ) \\
&= 4\pi^2 \big( ( 1 - {\rm Pr} (X\le-u) - {\rm Pr} (Y\le-v) + {\rm Pr}
(X\le-u,\, Y\le-v) ) \\
&\qquad- ( 1 - {\rm Pr} (X\le-u) ) ( 1 - {\rm Pr} (Y\le-v) ) \big) \\
&= 4\pi^2 ( {\rm Pr} (X\le-u,\, Y\le-v) - {\rm Pr}(X\le-u) {\rm
Pr}(Y\le-v) ).
\end{split}
\end{equation*}
By the same argument as in equation~(\ref{zazaza}) we see that
\begin{equation*}
\begin{split}
| \hat f(-u,-v) | &\le 4\pi^2 \min\{ {\rm Pr}(X\le-u), {\rm Pr}(Y\le-v) \} \\
&= 4\pi^2 \min\{ {\rm Pr}(X\ge u), {\rm Pr}(Y\ge v) \}
\end{split}
\end{equation*}
since $X$ and $Y$ are symmetric. We can therefore apply the
bounds~(\ref{thebounds}) to conclude that
\begin{equation*}
| \hat f(-u,-v) | \le 4\pi^2 \exp \big( {-0.04} (\max \{ |u|,|v| \} -3)^2
\big)
\end{equation*}
if either $|u|$ or $|v|$ exceeds 3. In particular, the function $\hat f$
decays (faster than) exponentially, as claimed above.

\subsection{Error due to truncating the infinite sums.}

From equation~(\ref{I8Error}) we have
\begin{equation}
\begin{split}
I_\qnn &= 4 \, \mathop{\sum \sum} \begin{Sb}  m, n \in \Z \\  m, n {\rm\ 
odd} \end{Sb}
\frac  { \hat \rho_\qnn ( m\epsilon/2, n\epsilon/2 )  
	- \hat \rho_\qnn ( m\epsilon/2, 0 ) \hat \rho_\qnn ( 0, n\epsilon/2 ) }
{mn}  \  + \  {\rm Error}_1 \\
&= \ 4 S_\qnn(\epsilon) \  + \  {\rm Error}_1
\end{split}
\label{omitaterm}
\end{equation}
where we have defined
\begin{equation}
S_\qnn(\epsilon) = \mathop{\sum \sum} \begin{Sb}  m, n \in \Z \\  m, n
{\rm\  odd} \end{Sb} \frac { \hat \rho_\qnn ( m\epsilon/2, n\epsilon/2 ) }
{mn}.
\label{Sepsilondef}
\end{equation}
(The term that has been omitted in the latter equality in
equation~(\ref{omitaterm}) equals zero, since $\hat \rho_\qnn (
m\epsilon/2, 0 ) \hat \rho_\qnn ( 0, n\epsilon/2 ) (mn)^{-1}$ is odd in
either variable separately due to the symmetry of the functions $\hat
\rho_\qnn ( m\epsilon/2, 0 )$ and $\hat \rho_\qnn ( 0, n\epsilon/2 )$
through the origin.) At this point we have accomplished the first step of
converting our integral $I_\qnn$ into a discrete sum, with a manageable
error; the next step is to truncate the ranges of summation so that the
resulting sum has only finitely many terms.

From the formula~(\ref{rhohatqnnformula}) for
$\hat\rho_\qnn$, the definition~(\ref{Sepsilondef}) becomes
\begin{equation}
\begin{split}
S_\qnn(\epsilon) &= \mathop{\sum \sum} \begin{Sb}  m, n \in \Z \\  m, n {\rm\ 
odd} \end{Sb} \frac
{F(m\epsilon,\chi_{-8}) F((n-m)\epsilon,\chi_{-4})
F(-n\epsilon,\chi_8)} {mn} \\
&= S_\qnn(\epsilon,C) + {\rm Error}_2,
\end{split}
\label{SqaaeCdef}
\end{equation}
where we have defined the truncated series
\begin{equation}
S_\qnn(\epsilon,C) = \mathop{\sumprime \sumprime}_{|m|, |n| \leq
2C/\epsilon} \frac {F(m\epsilon,\chi_{-8}) F((n-m)\epsilon,\chi_{-4})
F(-n\epsilon,\chi_8)} {mn}
\label{SCSS}
\end{equation}
(the primes indicating that the sums are taken over only odd values of
$m$ and $n$) and the error due to truncation
\begin{equation*}
{\rm Error}_2 = \mathop{\sumprime \sumprime}_{\max\{|m|, |n|\} > 2C/\epsilon}
\frac {F(m\epsilon,\chi_{-8}) F((n-m)\epsilon,\chi_{-4}) F(-n\epsilon,\chi_8)}
{mn}.
\end{equation*}
We rewrite this as
\begin{multline}
{\rm Error}_2 \  =  2 \bigg\{
  \sumprime_{m > 2C/\epsilon}  {\sumprime_{n=-m}^{m/2}}
+ \sumprime_{m > 2C/\epsilon}  {\sumprime_{n=m/2}^{m}} \\
+ \sumprime_{n > 2C/\epsilon} {\sumprime_{m=-n}^{n/2}}
+ \sumprime_{n > 2C/\epsilon} {\sumprime_{m=n/2}^{n}} \bigg\}
\frac {F(m\epsilon,\chi_{-8}) F((n-m)\epsilon,\chi_{-4})
F(-n\epsilon,\chi_8)} {mn},
\label{SCSSError}
\end{multline}
where the factor of 2 comes from grouping together the terms
corresponding to $(m,n)$ and $(-m,-n)$ by the symmetry of the summand
through the origin.

To bound Error$_2$, we will certainly need explicit estimates for the
functions $F(x,\chi)$ on the real axis. We recall the upper
bound~(\ref{FboundedwithJ}),
\begin{equation*}
|F(x,\chi)| \le (\pi|x|)^{-J/2} \prod_{j=1}^J \left ( \frac{1}{4} +
\gamma_j^2 \right )^{1/4},
\end{equation*}
where $J$ is any positive integer and $0<\gamma_1<\gamma_2<\dotsb$ are
the imaginary parts of the nontrivial zeros of $L(s,\chi)$. Any
particular choice of $J$ gives an upper bound of the form
\begin{equation}
| F(x,\chi) | \leq d(\chi) |x|^{-e(\chi)}  \label{Fxchibound}
\end{equation}
for some positive constants $d(\chi)$ and $e(\chi)$. For any fixed $x$
the optimal choice of $J$ is the largest integer such that $(\pi
x)^2>\frac14+\gamma_J^2$; but for our present purposes, we obtain
sufficiently good results that are easy to apply uniformly in $x$ by
choosing $J$ so that $\gamma_J$ is just less than~30. Table~\ref{dandeTable}
\begin{table}[bth]
\begin{tabular}{||c|c|c|c||}
\hline
\hline
$\chi$ & $J$ & $d(\chi)$ & $e(\chi)$ \\
\hline
\hline
$\chi_{-8}$ & 56 & $1.3\times10^{32}$ & 28 \\
\hline
$\chi_8$ & 56 & $2.1\times10^{32}$ & 28 \\
\hline
$\chi_{-4}$ & 46 & $8.5\times10^{26}$ & 23 \\
\hline
$\chi_{-3}$ & 42 & $7.5\times10^{24}$ & 21 \\
\hline
$\chi_{12}$ & 62 & $3.0\times10^{35}$ & 31 \\
\hline
\hline
\end{tabular}
\medskip
\caption{Allowable constants in the bound~(\ref{Fxchibound}) for
$|F(x,\chi)|$}
\label{dandeTable}
\end{table}
lists, for each of the five characters $\chi$ relevant to the
densities mod~8 and mod~12, the values of $J$ chosen and the resulting
values of $d(\chi)$ and $e(\chi)$ in the bound~(\ref{Fxchibound}),
which we computed from the lists of zeros of the $L(s,\chi)$ supplied
to us by R.~Rumely.

Since $|F(x,\chi)|$ is also bounded by 1 on the real axis, we can
estimate the first double sum in equation~(\ref{SCSSError}) by
\begin{equation}
\begin{split}
\bigg | \sumprime_{m > 2C/\epsilon} \sumprime_{n=-m}^{m/2} & \frac{
F(m\epsilon,\chi_{-8}) F((n-m)\epsilon,\chi_{-4}) F(-n\epsilon,\chi_8) }
{mn} \bigg | \\
&\leq \sum_{m > 2C/\epsilon} \sum_{n=-\infty}^{m/2} \left | \frac{
F(m\epsilon,\chi_{-8}) F((n-m)\epsilon,\chi_{-4}) } {m} \right | \\
&\leq   d(\chi_{-8}) d(\chi_{-4}) \epsilon^{-e(\chi_{-8})-e(\chi_{-4})}
\sum_{m > 2C/\epsilon} \sum_{n=-\infty}^{m/2}
m^{-e(\chi_{-8})-1}  (m-n)^{-e(\chi_{-4})}
\end{split}
\label{picktherightFs}
\end{equation}
using the bound~(\ref{Fxchibound}) for $\chi_{-8}$ and
$\chi_{-4}$.

Now we claim that
\begin{multline}
\sum_{m > M} \sum_{n=-\infty}^{m/2} m^{-\alpha} (m-n)^{-\beta} =
\sum_{m > M} \sum_{n=m/2}^\infty m^{-\alpha} n^{-\beta} \\
< 2^{\beta-1} M^{1-\alpha-\beta} \left ( \frac{2}{\alpha+\beta-1} +
\frac{M}{(\alpha + \beta -2)(\beta-1)} \right )
\label{jointbound}
\end{multline}
for any real numbers $\alpha,\beta>1$. The equality is clear upon
making the change of variables $n\mapsto m-n$, while the inequality
follows from the elementary argument%
\goodbreak%
\begin{equation*}
\begin{split}
\sum_{m > M} \sum_{n=m/2}^\infty m^{-\alpha} n^{-\beta} &< \sum_{m >
M} m^{-\alpha} \left ( \big\lceil \frac{m}{2} \big\rceil^{-\beta} +
\int_{ \lceil m/2 \rceil }^\infty t^{-\beta} \, dt \right ) \\
&\leq \sum_{m > M} m^{-\alpha} \left ( \left ( \frac{m}{2} \right
)^{-\beta} + \frac{1}{\beta-1} \left ( \frac{m}{2} \right )^{1-\beta}
\right ) \\
&< 2^\beta \int_M^\infty t^{-\alpha -\beta} \, dt +
\frac{2^{\beta-1}}{\beta-1} \int_M^\infty t^{1-\alpha -\beta} \, dt \\
&= \frac{2^\beta}{\alpha + \beta -1} M^{1-\alpha-\beta} +
\frac{2^{\beta-1}}{(\beta-1)(\alpha + \beta -2)} M^{2-\alpha-\beta},
\end{split}
\end{equation*}%
\goodbreak\noindent%
this last expression being equivalent to the right-hand side
of~(\ref{jointbound}).

Applying the upper bound~(\ref{jointbound}), with
$\alpha=e(\chi_{-8})+1$ and $\beta=e(\chi_{-4})$, to
equation~(\ref{picktherightFs}) gives
\begin{equation}
\begin{split}
\bigg | \sumprime_{m > 2C/\epsilon} \sumprime_{n=-m}^{m/2} & \frac{
F(m\epsilon,\chi_{-8}) F((n-m)\epsilon,\chi_{-4}) F(-n\epsilon,\chi_8) }
{mn} \bigg | \\
&< d(\chi_{-8}) d(\chi_{-4}) \epsilon^{-e(\chi_{-8})-e(\chi_{-4})} 2^{e(\chi_{-4})-1}
\left ( \frac{2C}{\epsilon} \right )^{-e(\chi_{-8})-e(\chi_{-4})} \\
&\qquad\times \left ( \frac{2}{e(\chi_{-8})+e(\chi_{-4})} +
\frac{2C/\epsilon}{(e(\chi_{-8})+e(\chi_{-4})-1)(e(\chi_{-4})-1)} \right ) \\
\vphantom{(^\big(} &= d(\chi_{-8}) d(\chi_{-4}) 2^{-e(\chi_{-8})}
C^{-e(\chi_{-8})-e(\chi_{-4})} \\
&\qquad\times \left ( \frac{1}{e(\chi_{-8})+e(\chi_{-4})} +
\frac{C/\epsilon}{(e(\chi_{-8})+e(\chi_{-4})-1)(e(\chi_{-4})-1)} \right ) \, .
\end{split}
\label{biglongeq}
\end{equation}
Substituting the appropriate values from Table~\ref{dandeTable}, we find
that this last expression is less than $1.85\times10^{-7}$ when
$\epsilon=1/20$ and $C=15$.

The second double sum in equation~(\ref{SCSSError}) may be similarly
bounded as
\begin{equation*}
\begin{split}
\bigg| \sumprime_{m > 2C/\epsilon} \sumprime_{n=-m/2}^{m} &
{F(m\epsilon,\chi_{-8}) F((n-m)\epsilon,\chi_{-4})
F(-n\epsilon,\chi_8) \over mn} \bigg| \\
&\leq \sum_{m > 2C/\epsilon} \sum_{n=m/2}^{\infty} \left | \frac{
F(m\epsilon,\chi_{-8}) F(-n\epsilon,\chi_8) } {m} \right | \\
&\leq d(\chi_{-8}) d(\chi_8) \epsilon^{-e(\chi_{-8})-e(\chi_8)}
\sum_{m > 2C/\epsilon} \sum_{n=m/2}^{\infty}
m^{-e(\chi_{-8})-1}n^{-e(\chi_8)}.
\end{split}
\end{equation*}
Applying~(\ref{Fxchibound}) to bound this last expression yields just
the right-hand side of equation~(\ref{biglongeq}) except with
$\chi_{-4}$ replaced with $\chi_8$; upon substituting
values from Table~\ref{dandeTable} we find that this expression is
also less than $1.85\times10^{-7}$ when $\epsilon=1/20$ and
$C=15$. The third and fourth double sums in~(\ref{SCSSError}) are
treated the same way, and so we conclude from
equation~(\ref{SCSSError}) that
\begin{equation}
|{\rm Error}_2| < 8(1.85\times10^{-7}) < 1.5\times10^{-6}
\label{Error2bound}
\end{equation}
when $\epsilon=1/20$ and $C=15$.

\subsection{Error due to approximating $F(z,\chi)$ by $F_T (z,\chi)$.}
We have accomplished the second step of approximating the infinite sum
$S_\qnn(\epsilon)$ by the finite sum $S_\qnn(\epsilon,C)$; however, this
latter sum is still unsuitable for computation, since it involves the
functions $F(z,\chi)$ which are infinite products. The last step is to
replace the functions $F(z,\chi)$ by their truncated counterparts
$F_T(z,\chi)$ defined in equation~(\ref{FTzchidefn}).

Recall the definition~(\ref{b1def}) of $b_1$,
\begin{equation*}
b_1=b_1(T,\chi) = - \sum_{\gamma \geq T} \frac {1} {\frac{1}{4} +
\gamma^2},
\end{equation*}
and put
\begin{equation}
\Delta_T (z,\chi) = \frac {\prod_{\gamma > T} J_0(\alpha_\gamma z) } { 1
+ b_1 z^2 } - 1.
\label{Delta}
\end{equation}
From the definitions~(\ref{Fzchidefn}) and~(\ref{FTzchidefn}) of $F$ and
$F_T$ we see that
\begin{equation*}
F (z,\chi) = F_T (z,\chi) (  1 + \Delta_T (z,\chi) ).
\end{equation*}
Making this substitution in equation~(\ref{SCSS}) for $\chi_{-8}$,
$\chi_{-4}$, and $\chi_8$, we then obtain
\begin{equation}
S_\qnn(\epsilon,C) = S_\qnn(\epsilon,C,T) + {\rm Error}_3,
\label{SqaaeCTdef}
\end{equation}
where $S_\qnn(\epsilon,C,T)$ is as defined in equation~(\ref{SCT}) and
\begin{multline}
{\rm Error}_3 = \mathop{\sumprime \sumprime}_{|m|, |n| \leq
2C/\epsilon} \frac {F_T(m\epsilon,\chi_{-8})
F_T((n-m)\epsilon,\chi_{-4}) F_T(-n\epsilon,\chi_8)} {mn} \\
{}\times \big( ( 1 + \Delta_T(m\epsilon,\chi_{-8}) ) ( 1 +
\Delta_T((n-m)\epsilon,\chi_{-4}) ) ( 1 + \Delta_T(-n\epsilon,\chi_8)
) - 1 \big).
\label{Err3expand}
\end{multline}

Regarding the size of the function $\Delta_T$, Rubinstein and Sarnak
\cite[Section 4.3]{RS} established the inequality
\begin{equation*}
\bigg| \bigg( \prod_{\gamma > T} J_0(\alpha_\gamma x) \bigg) - (1+b_1 x^2)
\bigg| \le {b_1^2 x^4 \over 2(1-|b_1| x^2)}
\end{equation*}
for real numbers $x$ satisfying $|b_1| x^2 < 1$. From the
definition~(\ref{Delta}) of $\Delta_T$ this immediately yields
\begin{equation}
| \Delta_T(x,\chi) | = {|\prod_{\gamma > T} J_0(\alpha_\gamma x) -
(1+b_1 x^2)| \over |1+b_1 x^2|} < {b_1^2 x^4 \over 2(1-|b_1| x^2)^2}
\qquad\hbox{if }|b_1| x^2 < 1.
\label{DeltaTbound}
\end{equation}

The quantities $b_1$ can be computed if we know all the zeros of
$L(s,\chi)$ up to height $T$, since
\begin{equation*}
b_1 = \sum_{0<\gamma<T} \frac {1} {\frac{1}{4} + \gamma^2} -
\sum_{\gamma>0} \frac {1} {\frac{1}{4} + \gamma^2}
\end{equation*}
and we have the formula (see Davenport \cite[p.~83]{davenport})
\begin{equation}
\sum_{\gamma>0} \frac {1} {\frac{1}{4} + \gamma^2} = \frac12 \sum_{\gamma}
\frac {1} {\frac{1}{4} + \gamma^2} = \frac12\log\frac q\pi -
{\gamma_0\over2} -(1+\chi(-1)){\log2\over2} + {L'(1,\chi)\over L(1,\chi)}
\label{Davenportformula}
\end{equation}
for a real primitive character $\chi\bmod q$, where
$\gamma_0=0.577215\dots$\ is Euler's constant.  The values $L(1,\chi)$
can be calculated in closed form by classical formulas (again
see~\cite{davenport}), while the values $L'(1,\chi)$ can be calculated
in closed form using a formula of Selberg and Chowla~\cite{SC} for the
odd characters and a formula of Deninger~\cite{deninger} for the even
characters. The former formula expresses $L'(1,\chi)$ in terms of the
logarithm of the $\Gamma$-function, while the latter expresses
$L'(1,\chi)$ in terms of a function $R(x)$ defined as
\begin{equation*}
R(x)= \big( {\partial^2\zeta(s,x)\over\partial s^2} \big) \big|_{s=0};
\end{equation*}
here $\zeta(s,x)$ is the Hurwitz zeta function, defined when $x>-1$ by
$\zeta(s,x) = \sum_{n=1}^\infty (n+x)^{-s}$ for $\Re s>1$ and by
meromorphic continuation to the complex $s$-plane.

The Mathematica software package
is capable of calculating both $\log\Gamma(x)$ and $R(x)$ to
arbitrary precision, and thus by the formula~(\ref{Davenportformula}) the
sums $\sum_{\gamma>0} 1/(\frac14 + \gamma^2)$ can also be so
calculated. Table~\ref{L1chitable} contains the results of such
calculations for the five characters relevant to the densities mod 8 and
mod 12.
\begin{table}[tb]
\begin{tabular}{||c|c|c|c||}
\hline
\hline
$\chi$ & $L(1,\chi)$ & $L'(1,\chi)$ & $\sum_{\gamma>0}{1\over1/4+\gamma^2}$ \\
\hline
\hline
$\chi_{-8}$ & ${\pi\over2\sqrt2}$ & ${\pi\over2\sqrt2}(\gamma_0+\log2\pi+\log{\Gamma(5/8)\Gamma(7/8)\over\Gamma(1/8)\Gamma(3/8)})$ & $0.158037$ \\
\hline
$\chi_8$ & ${\log(1+\sqrt2)\over\sqrt2}$ & ${1\over2\sqrt2}(\gamma_0+\log2\pi+R({1\over8})-R({3\over8})-R({5\over8})+R({7\over8}))$ & $0.117716$ \\
\hline
$\chi_{-4}$ & $\frac\pi4$ & $\frac\pi4(\gamma_0+\log2\pi+2\log{\Gamma(3/4)\over\Gamma(1/4)})$ & $0.077784$ \\
\hline
$\chi_{-3}$ & ${\pi\over3\sqrt3}$ &
${\pi\over3\sqrt3}(\gamma_0+\log2\pi+3\log{\Gamma(2/3)\over\Gamma(1/3)})$ &
$0.056615$ \\
\hline
$\chi_{12}$ & ${\log(2+\sqrt3)\over\sqrt3}$ & ${1\over2\sqrt3}(\gamma_0+\log2\pi+R({1\over12})-R({5\over12})-R({7\over12})+R({11\over12}))$ & $0.165083$ \\
\hline
\hline
\end{tabular}
\medskip
\caption{Values of $L(1,\chi)$, $L'(1,\chi)$, and
$\sum_{\gamma>0}{1\over1/4+\gamma^2}$}
\label{L1chitable}
\end{table}

For all five of these characters, when we choose $T=1$0,000 we find that
$|b_1|<0.000173$. The upper bound~(\ref{DeltaTbound}) can then be written
more simply as $|\Delta_T(x,\chi)|\le D(x)$ for $|x|<74$, where we have
defined
\begin{equation}
D(x) = \frac{1.5 \times 10^{-8}x^4}{(1-0.00018 x^2)^2}.
\label{Dzdef}
\end{equation}
Consequently, the definition~(\ref{Err3expand}) of Error$_3$ implies
the inequality 
\begin{multline}
\left | {\rm Error}_3 \right | \leq \mathop{\sumprime \sumprime}_{|m|,
|n| \leq 2C/\epsilon} \left | \frac {F_T(m\epsilon,\chi_{-8})
F_T((n-m)\epsilon,\chi_{-4}) F_T(-n\epsilon,\chi_8)} {mn} \right | \\
\times \big( ( 1 + D(m\epsilon) ) ( 1 + D((n-m)\epsilon) ) ( 1 +
D(-n\epsilon) ) - 1 \big).
\label{error3ugly}
\end{multline}
The quantity on the right-hand side of this inequality was computed at
the same time as the sum $S_\qnn(\epsilon,C,T)$ was computed,
and we obtained the bound
\begin{equation}
|{\rm Error}_3| < 5.5\times10^{-6}.  \label{Error3bound}
\end{equation}

\subsection{Conclusion.}

From the relationships~(\ref{Sepsilondef}),~(\ref{SqaaeCdef}),
and~(\ref{SqaaeCTdef}) among the various intermediate sums $S_\qnn$, we have
\begin{equation*}
I_\qnn = 4 S_\qnn(\epsilon,C,T) + {\rm Error}_1 + 4 {\rm Error}_2 + 4 {\rm
Error}_3.
\end{equation*}
Using this identity in equation~(\ref{the8357formula}) yields
\begin{equation*}
\delta_\qnn = \frac{1}{4} - \frac{1}{4\pi^2} I_\qnn = \frac14 -
{1\over4\pi^2} ( 4S_\qnn(\epsilon,C,T) + {\rm Error}_1 + 4{\rm Error}_2 +
4{\rm Error}_3 ),
\end{equation*}
whence it follows that
\begin{equation*}
\big| \delta_\qnn - \big( \frac14 -
{S_\qnn(\epsilon,C,T)\over\pi^2} \big) \big| \le {|{\rm
Error}_1|\over4\pi^2} + {|{\rm Error}_2|+|{\rm Error}_3|\over\pi^2}.
\end{equation*}
Thus by the inequalities~(\ref{Error1bound}),~(\ref{Error2bound}),
and~(\ref{Error3bound}), we conclude that
\begin{equation*}
\left | \delta_\qnn - \big( \frac14 - {S_\qnn(\epsilon,C,T)\over\pi^2} \big)
\right | <  8\times10^{-7}
\end{equation*}
when $\epsilon=1/20$, $C=15$, and $T=1$0,000.  Using these values for
$\epsilon$, $C$, and $T$, the sum $S_\qnn(\epsilon,C,T)$ was
calculated and found to equal $0.5645285\dots$, and therefore we have
rigorously that
\begin{equation*}
\delta_\qnn = 0.1928013 \pm 9\times10^{-7},
\end{equation*}
which is slightly stronger than the first assertion of Theorem 1.

The error analysis in Sections~3.1--3.3 can be repeated for each of
the densities in Theorem~1; the constants mentioned in the
error analysis have been chosen to apply to all of these densities.
Therefore, the densities calculated for Theorem 1 are all correct to within
the same margin $9\times10^{-7}$, which is enough to establish the theorem.

\subsection{Appendix: Probability bounds.}

In this section we establish the bounds~(\ref{thebounds}) for
${\rm Pr} ( X \geq u )$ and ${\rm Pr} ( Y \geq u )$
which were used for the computations in Section 3.1.

To do so, we first recall from Section 2.1 the explicit form of the random
variables having the distribution $\mu_\qaa$. Specializing the
representation~(\ref{bXchi}) to the case $q=8$ and $\{a_1,a_2,a_3\} =
\{3,5,7\}$, we find that $\mu_\qnn$ is the distribution of the random
$\R^3$-vector
\begin{equation*}
(1,1,1)  +  X(\chi_{-8}) ( 1,-1,-1 ) +  X(\chi_{-4}) ( -1,1,-1  )
+  X(\chi_8) ( -1,-1,1 ).
\end{equation*}
Next, recalling the changes of variables~(\ref{mu2nu})
and~(\ref{nu2rho}) that took us from $\mu$ to $\nu$ and then to
$\rho$, we observe that $\rho_\qnn$ is the distribution of of the
random $\R^2$-vector
\begin{equation}
X(\chi_{-8}) ( 2,0 ) +  X(\chi_{-4}) ( -2,2 ) +  X(\chi_8) ( 0,-2 ) \, .
\label{Xrho}
\end{equation}
Now define the two real-valued random variables
\begin{equation}
\begin{split}
X &\ =\  2 \sum \begin{Sb}\gamma>0 \\ L(1/2+i\gamma,\chi_{-8})=0 \end{Sb}
\alpha_\gamma \sin ( 2\pi U_\gamma )\ - \ 2 \sum \begin{Sb}\gamma>0 \\
L(1/2+i\gamma,\chi_{-4})=0 \end{Sb} \alpha_\gamma \sin ( 2\pi U_\gamma ),
\\
Y &\ =\  2 \sum \begin{Sb}\gamma>0 \\ L(1/2+i\gamma,\chi_{-4})=0 \end{Sb}
\alpha_\gamma \sin ( 2\pi U_\gamma ) \ - \ 2 \sum \begin{Sb}\gamma>0 \\
L(1/2+i\gamma,\chi_8)=0 \end{Sb} \alpha_\gamma \sin ( 2\pi U_\gamma ).
\end{split}
\label{XandYdef}
\end{equation}
We see from the definition~(\ref{Xchi}) of the $X(\chi)$ that the
random vector $(X,Y)$ equals the random vector~(\ref{Xrho}).

The following lemma gives information about the tails of
random variables of this type.
\medskip

\begin{lemma}
Let $r_1\ge r_2\ge\dotsb$ be a sequence of positive real numbers
such that $\sum_{k=1}^\infty r_k =\infty$ but
$\sum_{k=1}^\infty r_k^2 = R < \infty$. Let $U_1$, $U_2$, \dots be
independent random variables uniformly distributed on $[0,1]$, and define
the random variable
\begin{equation*}
W = \sum_{k=1}^\infty r_k \sin ( 2\pi U_k ).
\end{equation*}
Then for any real number $w\ge2r_1$,
\begin{equation*}
{\rm Pr}(W\ge w) \le \exp \big( {-3(w-2r_1)^2 \over 16R}
\big).
\end{equation*}
\label{montgomerylemma}
\end{lemma}

\noindent{\bf Proof.}
Theorem 1 of Montgomery \cite[Section 3]{montgomery}
states that under the assumptions of this lemma, we have
\begin{equation}
{\rm Pr}  \bigg( W \geq 2 \sum_{k=1}^K r_k \bigg)  \leq
\exp \bigg( {-\frac34}  \bigg( \sum_{k=1}^K r_k \bigg)^2 \big/
\sum_{k>K} r_k^2 \bigg)
\label{monty}
\end{equation}
for any integer $K \geq 1$. Since the $r_k$ are decreasing and
$\sum_{k=1}^\infty r_k = \infty$, it is clear that for any $w \geq 2r_1$
there exists a $K\ge1$ such that
\begin{equation*}
{w\over2} - r_1 \leq \sum_{k=1}^K r_k \leq {w\over2}.
\end{equation*}
With this choice of $K$, the inequality~(\ref{monty}) simplifies to
\begin{equation*}
{\rm Pr} (W\ge w) \le \exp \bigg( {-\frac34} \big( \frac w2 - r_1 \big)^2
\big/ \sum_{k=K}^\infty r_k^2 \bigg) \le \exp \big(
\frac{-3(w-2r_1)^2}{16R} \big)
\end{equation*}
which is the statement of the lemma. \qed
\medskip

We now apply this lemma to the random variables $X$ and $Y$ defined in
equation~(\ref{XandYdef}). (Note that because each variable $U_\gamma$
is uniformly distributed on $[0,1]$, we may replace each $U_\gamma$ in
the second sums on each line with $U_\gamma+1/2$; this has the effect
of changing the subtraction signs in the equations~(\ref{XandYdef}) to
addition signs, thus rendering $X$ and $Y$ into the form to which
Lemma 3.1 applies.)  For the variable $X$, the sequence corresponding
to $r_k$ is
\begin{equation*}
\{2\alpha_\gamma\colon L(1/2+i\gamma,\chi_{-8})=0,\, \gamma>0\} \cup
\{2\alpha_\gamma\colon L(1/2+i\gamma,\chi_{-4})=0,\, \gamma>0\}.
\end{equation*}
For this sequence, the largest element $r_1$ is less than $1.5$, and the
sum $R$ of the squares of the elements does not exceed $4.5$. Therefore,
applying Lemma~\ref{montgomerylemma}, we find that
\begin{equation*}
{\rm Pr} ( X \geq u ) \le \exp ({-0.04} (u-3)^2)
\end{equation*}
for any $u\ge3$. Similarly, $Y$ can be shown to satisfy the same
estimate, which establishes the upper bounds~(\ref{thebounds}). In fact, the
constants mentioned above will work for every pair of characters that
arises in the computations of $\rho_{8;a_1,a_2,a_3}$, where
$\{a_1,a_2,a_3\}$ is a permutation of $\{3,5,7\}$, and in
$\rho_{12;a_1,a_2,a_3}$, where $\{a_1,a_2,a_3\}$ is a permutation of
$\{5,7,11\}$.

\bigskip
\section{Computational Results.}
\setcounter{equation}{0}
\setcounter{table}{0}

\newcommand{\up}[1]{\smash{\raise7pt\hbox{#1}}}
\renewcommand{\center}[1]{\multicolumn{1}{||c|}{#1}}

The mathematical and numerical computations described in this paper
were implemented on an SGI Challenge computer using the Mathematica
software package, which has the capability to perform computations to
arbitrary and verifiable precision (see Wolfram~\cite{mathematica}). A
typical quantity to be calculated is the expression
$S_\qnn(\epsilon,C,T)$ defined in equation~(\ref{SCT}), which depends
on the functions $F_T(z,\chi)$ defined in equation~(\ref{FTzchidefn}).
To compute these functions we needed, for the Dirichlet $L$-functions
corresponding to characters to the moduli $q\le12$, lists of the zeros
whose imaginary parts are bounded by $T=1$0,000. These lists of
imaginary parts of zeros (accurate to twelve decimal places) were
kindly supplied to us by R.~Rumely (see~\cite{rumely}). For the
estimation of Error${}_3$ in Section 3.3 it was also necessary to
compute quantities typified by the right-hand side of
equation~(\ref{error3ugly}), which is no harder than computing
$S_\qnn(\epsilon,C,T)$ itself.

In addition to the results reported in Theorem 1, a number of further
computations were carried out involving certain cases with $q\le12$
and $r\le4$. In these additional results, which are presented below,
we report only the numbers of decimal places in which we have some
degree of confidence; specifically, we expect the entries to be
correct to within one or two units in the last decimal place reported.

Table~\ref{table2way345}
\begin{table}[b]
\begin{tabular}{||c|l|c||}
\hline \hline
$q$ & \multicolumn{1}{|c|}{$a_1a_2$} & $\delta_{q; a_1,a_2}$ \\ \hline \hline
  & NS: 21 & .9990633 \\
\up3 & SN: 12 & .0009367 \\
\hline
\hline
  & NS: 31 & .9959280 \\
\up4 & SN: 13 & .0040720 \\
\hline
\hline
 & NS: 21,24,31,34 & .952140 \\ \cline{2-3}
 & NN: 23,32 & \\
\up5& SS: 14,41 & \up{$1/2$} \\ \cline{2-3}
 & SN: 12,13,42,43 & .047860 \\
\hline
\hline
\end{tabular}
\medskip
\caption{Two-way races for the moduli $q=3,4,5$}
\label{table2way345}
\end{table}

shows the calculated densities $\delta_{q;a_1,a_2}$ for the two-way races
between $\pi(x;q,a_1)$ and $\pi(x;q,a_2)$, for the moduli
$q=3$, 4, and 5. For example, the first line of the table indicates that
$\delta_{3;2,1}=0.999063$ (rounded to seven decimal places). Throughout this
section we use the symbol $N$ to stand for any nonsquare mod $q$ and $S$ to
stand for any square mod $q$ (although distinct occurrences of $N$ or $S$ in
a single entry stand for distinct residues) to make the Chebyshev
biases more clearly evident where appropriate.

Of course, since $\varphi(3)=\varphi(4)=2$, the two-way races shown are the
only possible races for the moduli 3 and 4. The densities for these moduli
were calculated by Rubinstein and Sarnak and our calculations agree with
theirs to six decimal places. (Although they were only reported in~\cite{RS}
truncated to four decimal places, they had in fact been calculated to higher
accuracy.)

For the races modulo 5, it turns out that the densities $\delta_{q;a_1,a_2}$
depend only on whether or not $a_1$ and $a_2$ are squares mod 5, due to the
symmetry results given in Theorem 2. (In fact this is true for the races
between multiple residues mod 5 as well.) For instance, applying Theorem
2(b) with $a_1=2$, $a_2=1$, and $b=4$ shows that
$\delta_{5;2,1}=\delta_{5;3,4}$; then applying Theorem 2(a) to each of these
expressions shows further that $\delta_{5;2,1}=\delta_{5;3,1}$ and
$\delta_{5;3,4}=\delta_{5;2,4}$. Since the two nonsquares mod 5 are $\{2,3\}$
while the two squares are $\{1,4\}$, these equalities show that all four
densities represented by $\delta_{5;N,S}$ are equal, as indicated in
Table~\ref{table2way345}.

The fact that $\delta_{q;N,N}=\delta_{q;S,S}=1/2$ as shown on the
penultimate line of the table was proved by Rubinstein and Sarnak, and
it also follows from our Theorem 2(d). We calculated these densities
anyway, and the calculated answers differed from $1/2$ by at most
$10^{-16}$, which is the default machine precision for our Mathematica
calculations. This degree of accuracy is not unexpected in this
instance, as the integral in the formula~(\ref{genRS}) is identically
zero when $a_1$ and $a_2$ are both squares or both nonsquares mod~$q$.

Table~\ref{table3way5}
\begin{table}[bth]
\begin{tabular}{||l|c||}
\hline
\hline
\center{$a_1a_2a_3$} & $\delta_{5; a_1,a_2,a_3}$ \\
\hline
\hline
NNS: 231,234,321,324 & \\
NSS: 214,241,314,341 & \up{.45678} \\
\hline
NSN: 213,243,312,342 & \\
SNS: 124,134,421,431 & \up{.03859} \\
\hline
SNN: 123,132,423,432 & \\
SSN: 142,143,412,413 & \up{.00464} \\
\hline
\hline
\end{tabular}
\medskip
\caption{Three-way races modulo $q=5$}
\label{table3way5}
\end{table}

provides the calculated densities $\delta_{q;a_1,a_2,a_3}$ for the
three-way races modulo 5. Again, in this case the densities only
depend on whether $a_1$, $a_2$, and $a_3$ are squares mod 5, by the
symmetry results (a) and (b) of Theorem 2. In addition, each density
matches two different types of permutations: for instance, Theorem
2(e) with $a_1=2$, $a_2=3$, $a_3=1$, and $b=2$ asserts that
$\delta_{5;2,3,1}=\delta_{5;2,1,4}$ as indicated in the first entry of
the table.

As mentioned at the beginning of this section, we are confident from
numerical considerations that the numbers reported in Table~\ref{table3way5}
are accurate to the five decimal places given there, with a
possible error of one or two units in the fifth decimal place. Thus, for instance,
if we choose a particular triple of residues such as $\{1,2,3\}$ and add
up the densities from Table~\ref{table3way5} corresponding to the six
permutations of that triple, the result is $1.00002$. Moreover, the
three ordered triples $\{3,2,1\}$, $\{2,3,1\}$, and $\{2,1,3\}$ are
the three permutations in which 2 is ahead of 1, and so we have the
identity
\begin{equation}
\delta_{5;2,1} = \delta_{5;3,2,1} + \delta_{5;2,3,1} +
\delta_{5;2,1,3}
\label{deltasplit}
\end{equation}
(cf.~equation~(\ref{a3choices})). Table~\ref{table2way345} gives
$.952140$ for the left-hand side of this identity, while adding the
appropriate entries from Table~\ref{table3way5} gives
$.95215$ for the right-hand side.

There are two reasons why our calculations of the densities in
three-way races for moduli other than 8 and 12 are less accurate than
the full six-decimal-place accuracy proven in Theorem~1, both stemming from
the fact that there are complex-valued Dirichlet characters associated
with the other moduli. First, when we calculate the function
$F_T(z,\chi)$ we do so only on a discrete set of points, evenly
spaced at intervals of $\epsilon/2$. These points are the only ones
needed to evaluate the sum $S_\qnn(\epsilon,C,T)$, as we see from its
definition~(\ref{SCT}), but for the sums corresponding to other moduli we
need to know the value of $F_T(z,\chi)$ at irrational multiples
of~$\epsilon$. We estimated this value by interpolating linearly between the
two nearest values, and this estimation introduces an additional error into the
calculations.

Second, the zeros of $L$-functions corresponding to complex
characters are not symmetric with respect to the real axis, and so the
quantity $\sum_{\gamma>0} 1/(\frac14+\gamma^2)$, needed to compute
$b_1(T,\chi)$, cannot be evaluated in closed form. Since we can
evaluate $b_1(T,\chi)+b_1(T,\bar\chi)$ in closed form, we used half of
this quantity in place of both $b_1(T,\chi)$ and $b_1(T,\bar\chi)$;
this gives the correct first-order approximation to the tail of
$F(z,\chi)F(z,\bar\chi)$, but the absolute error in our
calculations can be somewhat higher as a result. For higher moduli,
the sheer number of characters will also play a role, as the product
of the $\phi(q)-1$ functions $F_T(z,\chi)$ required for the evaluation
of $\hat\rho_\qaa$ will gradually erode the accuracy of the calculated
number.

Since there are precisely four reduced residues modulo 5, it is
natural to look at the complete four-way race mod 5;
Table~\ref{table4way5}
\begin{table}[bth]
\begin{tabular}{||l|c||}
\hline \hline
\center{$a_1a_2a_3a_4$} & $\delta_{8; a_1,a_2,a_3}$ \\ \hline \hline
NNSS: 2314,2341,3214,3241 & $.21136$ \\ \hline
NSNS: 2134,2431,3124,3421 & $.02985$ \\ \hline
NSSN: 2143,2413,3142,3412 & \\
SNNS: 1234,1324,4231,4321 & \up{$.00424$} \\ \hline
SNSN: 1243,1342,4213,4312 & $.00028$ \\ \hline
SSNN: 1423,1432,4123,4132 & $.00007$ \\ \hline
\hline
\end{tabular}
\medskip
\caption{The full four-way race modulo $q=5$}
\label{table4way5}
\end{table}

shows the calculated densities for this four-way race. Here again,
the densities only depend on whether $a_1$, $a_2$, $a_3$,
and $a_4$ are squares mod 5, by the symmetry results from parts (a) and
(b) of Theorem 2, with the added symmetry in the third entry of the
table following from Theorem 2(e). Once again we can estimate the
accuracy of these densities by comparing the sum of all twenty-four
densities to 1, and also by comparing the values here to those in
Table~\ref{table3way5} using identities such as
\begin{equation*}
\delta_{5;1,2,3} = \delta_{5;4,1,2,3} + \delta_{5;1,4,2,3} +
\delta_{5;1,2,4,3} + \delta_{5;1,2,3,4}.
\end{equation*}
In all cases, these sums of densities from Table~\ref{table4way5} are
precise to within a few units in the fifth decimal place.

In the calculation of these four-way densities, the general formula given
in Theorem 4 involves a three-dimensional integral which must be computed
numerically. Performing this calculation with a reasonable degree of
accuracy lies at the limit of the computing capabilities of the method used for the
calculations in this paper; in particular, we found it necessary to reduce
the value of $C$ and increase the value of $\epsilon$ somewhat to
make the computations feasible.

Since the distribution of the primes into residue classes modulo 6 is
fully determined by their distribution mod 3, the next modulus of
interest is $q=7$. Table~\ref{table2way7}
\begin{table}[b]
\begin{tabular}{||c|c||}
\hline
\hline
${a_1a_2}$ & $\delta_{7; a_1,a_2}$ \\
\hline
\hline
31,32,51,54,62,64 & .874349 \\
\hline
34,52,61 & .845210 \\
\hline
12,14,21,24,41,42, & \\
\,35,36,53,56,63,65 & \up{$1/2$} \\
\hline
16,25,43 & .154790 \\
\hline
13,15,23,26,45,46 & .125651 \\
\hline
\hline
\end{tabular}
\medskip
\caption{Two-way races modulo $q=7$}
\label{table2way7}
\end{table}

shows the calculated densities $\delta_{7;a_1,a_2}$ for the two-way
races modulo 7. Here for the first time, we see that the density does
not depend merely on whether $a_1$ and $a_2$ are squares mod 7: the
squares mod 7 are $\{1,2,4\}$, and so each of the top two lines of the
table are densities of the form $\delta_{7;N,S}$, while each of the
bottom two lines are densities of the form $\delta_{7;S,N}$. In other
words, Chebyshev's bias is not the only factor causing
asymmetries in the Shanks--R\'enyi race games. (For a somewhat more
precise discussion of Chebyshev biases for $r$-tuples with $r\ge3$, see
the discussion of ``bias factors'' in Section 6.) The middle row of the
table again indicates the known fact that all densities of the form
$\delta_{7;N,N}$ and $\delta_{7;S,S}$ equal $1/2$.

Table~\ref{table3way7}
\begin{table}[bth]
\begin{tabular}{||c|c||}
\hline
\hline
\center{$a_1a_2a_3$} & $\delta_{7; a_1,a_2,a_3}$ \\
\hline
\hline
512; 314; 631; 651; 621; 324; 532; 562; 641; 542; 354; 364 & .4038 \\ 
521; 341; 361; 561; 612; 342; 352; 652; 614; 524; 534; 634 & .3678 \\ 
251; 431; 316; 516; 162; 432; 325; 625; 164; 254; 543; 643 & .1027 \\ 
152; 134; 613; 615; 261; 234; 523; 526; 461; 452; 345; 346 & .0736 \\ 
215; 413; 136; 156; 126; 423; 235; 265; 146; 245; 453; 463 & .0295 \\ 
125; 143; 163; 165; 216; 243; 253; 256; 416; 425; 435; 436 & .0226 \\ 
\hline
\hline
312,321; 351,531; 514,541; 362,632; 624,642; 564,654 & .3943 \\ 
132,231; 315,513; 154,451; 326,623; 264,462; 546,645 & .0857 \\ 
123,213; 135,153; 145,415; 236,263; 246,426; 456,465 & .0200 \\ 
\hline
\hline
124,142,214,241,412,421; 356,365,536,563,635,653 & $1/6$ \\
\hline
\hline
\end{tabular}
\medskip
\caption{Three-way races modulo $q=7$}
\label{table3way7}
\end{table}

gives the calculated densities for the three-way races modulo
7. Because the number of different values for the densities is larger
than in the previous cases, we have not organized them strictly by
decreasing size, but rather we have grouped together the values
corresponding to {\it isomorphic race games\/}. We say that two $r$-tuples
$\{a_1,\dots,a_r\}$ and $\{b_1,\dots,b_r\}$ of reduced residue classes mod
$q$ have isomorphic race games if there exists a bijection $\tau$ from the
set $\{1,\dots,n\}$ to itself such that each residue $a_j$ acts exactly
like the corresponding residue $b_{\tau(j)}$, i.e., if
\begin{equation*}
\delta_{q;a_{\sigma(1)},\dots,a_{\sigma(r)}} =
\delta_{q;b_{\tau(\sigma(1))},\dots,b_{\tau(\sigma(r))}}
\end{equation*}
for every permutation $\sigma$ of $\{1,\dots,n\}$.

For instance, Theorem 2(a) tells us that
$\delta_{7;1,2,5}=\delta_{7;1,3,4}$ and similarly for the corresponding
permutations of $\{1,2,5\}$ and $\{1,3,4\}$. Therefore the bijection
$\tau\colon\{1,2,5\}\to\{1,3,4\}$ given by $\tau(a)\equiv a^{-1}\pmod7$
shows that these triples have isomorphic race games. Table~\ref{table3way7}
shows that there are ten triples whose race games are in the isomorphism
class determined by $\{1,2,5\}$; the six densities for the race games in this
class are all distinct. In addition, there are five triples in the isomorphism
class of $\{1,2,3\}$; the race games in this class have only three distinct
densities due to an additional symmetry generated by Theorem 2(a). Finally,
the two special triples \{S,S,S\}${}=\{1,2,4\}$ and \{N,N,N\}${}=\{3,5,6\}$
each give completely symmetric race games; this is the smallest modulus to
which parts (d) and (e) of Theorem 2 can be applied, since three distinct
squares or nonsquares are needed. The complete symmetry for these two race
games was also proven by Rubinstein and Sarnak. We remark that our
computations of these densities yielded $1/6$ to five decimal places. We
did not proceed further with computations modulo 7, since there is no
natural four-way race and races with five or more residues are beyond the
present capabilities of our computing set-up.

Table~\ref{table2way8}
\begin{table}[bth]
\begin{tabular}{||c|c||}
\hline \hline
${a_1a_2}$ & $\delta_{8; a_1,a_2}$ \\ 
\hline
\hline
31 & .9995688 \\
13 & .0004312 \\
\hline
\hline
51 & .9973946 \\
15 & .0026054 \\
\hline
\hline
71 & .9989378 \\
17 & .0010622 \\
\hline
\hline
35,37,53, & \\
57,73,75 & \up{$1/2$} \\
\hline 
\hline
\end{tabular}
\medskip
\caption{Two-way races modulo $q=8$}
\label{table2way8}
\end{table}

shows the calculated densities for the two-way races modulo 8. Because
only one fourth of the residues mod 8 are squares (i.e., $c(8,1)=3$),
in contrast to the lower moduli, there are fewer symmetries among the
densities. (This is somewhat counterintuitive, since the
multiplicative group modulo 8 is highly symmetric.) This higher value
of $c(8,1)$ also causes a larger bias towards nonsquares, as can be
seen by the fact that the values in Table~\ref{table2way8} are more
extreme than those in Tables~\ref{table2way345} and~\ref{table2way7}.

Table~\ref{table3way8}
\begin{table}[bth]
\begin{tabular}{||c|c||c|c||c|c||c|c||}
\hline \hline
$a_1a_2a_3$ & $\delta_{8; a_1,a_2,a_3}$ &
$a_1a_2a_3$ & $\delta_{8; a_1,a_2,a_3}$ &
$a_1a_2a_3$ & $\delta_{8; a_1,a_2,a_3}$ &
$a_1a_2a_3$ & $\delta_{8; a_1,a_2,a_3}$ \\
\hline
\hline
$531$ & $.4996015$ & $731$ & $.4995765$ & $571$ & $.4990135$ & & \\
$351$ & $.4974123$ & $371$ & $.4989440$ & $751$ & $.4974474$ & \up{357,753} & \up{$.1928013$} \\
\cline{7-8} 
$315$ & $.0025550$ & $317$ & $.0010483$ & $715$ & $.0024769$ & & \\
$513$ & $.0003808$ & $713$ & $.0004173$ & $517$ & $.0009337$ & \up{375,573} & \up{$.1664263$} \\
\cline{7-8}
$135$ & $.0000327$ & $137$ & $.0000077$ & $175$ & $.0000757$ & & \\
$153$ & $.0000177$ & $173$ & $.0000062$ & $157$ & $.0000528$ & \up{735,537} & \up{$.1407724$} \\
\hline
\hline 
\end{tabular}
\medskip
\caption{The four three-way races modulo $q=8$}
\label{table3way8}
\end{table}

shows the calculated densities for the three-way races modulo 8, including
the values for $\delta_{8;N,N,N}$ highlighted in Theorem 1. Since all of
the characters mod 8 are real, the additional sources of computational
error mentioned in the discussion of Table~\ref{table3way5} are not present
here, and so we feel justified in reporting these figures to seven decimal
places; in fact note that the sums of the
appropriate three-way densities sum to the two-way densities in
Table~\ref{table2way8} in a manner analogous to
equation~(\ref{deltasplit}), with the sums all agreeing to within one
or two units in the seventh decimal place.

As with the modulus 5, it is natural to look at the complete four-way race
modulo 8;
Table~\ref{table4way8}
\begin{table}[b]
\begin{tabular}{||c|c|c|c|c|c|c|c||}
\hline \hline
$a_1a_2a_3a_4$ & $\delta_{8; a_1,a_2,a_3,a_4}$ &
$a_1a_2a_3a_4$ & $\delta_{8; a_1,a_2,a_3,a_4}$ &
$a_1a_2a_3a_4$ & $\delta_{8; a_1,a_2,a_3,a_4}$ &
$a_1a_2a_3a_4$ & $\delta_{8; a_1,a_2,a_3,a_4}$ \\ \hline \hline
$1357$ & $.0000014$ & $3157$ & $.0000500$ & $5137$ & $.0000027$ & $7135$ & $.0000261$ \\
$1375$ & $.0000029$ & $3175$ & $.0000696$ & $5173$ & $.0000023$ & $7153$ & $.0000154$ \\
$1537$ & $.0000007$ & $3517$ & $.0007972$ & $5317$ & $.0001315$ & $7315$ & $.0008983$ \\
$1573$ & $.0000006$ & $3571$ & $.1919526$ & $5371$ & $.1406374$ & $7351$ & $.1398456$ \\
$1735$ & $.0000023$ & $3715$ & $.0015371$ & $5713$ & $.0000848$ & $7513$ & $.0002910$ \\
$1753$ & $.0000009$ & $3751$ & $.1648170$ & $5731$ & $.1663386$ & $7531$ & $.1924939$ \\
\hline
\hline 
\end{tabular}
\medskip
\caption{The full four-way race modulo $q=8$}
\label{table4way8}
\end{table}

shows the calculated densities for this four-way race, listed in the
lexicographical ordering on the permutations of $\{1,3,5,7\}$. Despite the
need to use slightly cruder values of $C$ and $\epsilon$ in the
calculations of the three-dimensional integrals arising in the formulas for
these densities, the sum of all 24 densities and numerical checks against
Table~\ref{table3way8} suggest that these densities are also accurate to
within one or two units in the seventh decimal place.

Tables~\ref{table2way9} and~\ref{table3way9}
\begin{table}[bth]
\begin{tabular}{||c|c||}
\hline
\hline
\center{$a_1a_2$} & $\delta_{9; a_1,a_2}$ \\
\hline
\hline
21,24,51,57,84,87 & .881584 \\
\hline
27,54,81 & .864230 \\
\hline
14,17,41,47,71,74, & \\
25,28,52,58,82,85 & \up{$1/2$} \\
\hline
72,45,18 & .135770 \\
\hline
12,15,42,48,75,78 & .118416 \\
\hline
\hline
\end{tabular}
\medskip
\caption{Two-way races modulo $q=9$}
\label{table2way9}
\end{table}

show the calculated densities for the two-way and three-way races
modulo 9. Since the multiplicative group mod 9 is isomorphic to the
multiplicative group mod 7 (both are cyclic of order 6), the various
symmetries present in Tables~\ref{table2way9} and~\ref{table3way9}
mirror those found in Tables~\ref{table2way7} and~\ref{table3way7},
with the squares mod 9 being $\{1,4,7\}$.
\begin{table}[t]
\begin{tabular}{||c|c||}
\hline
\hline
\center{$a_1a_2a_3$} & $\delta_{9; a_1,a_2,a_3}$ \\
\hline
\hline
514; 217; 821; 851; 841; 247; 524; 584; 871; 574; 257; 287 & .4010 \\
541; 271; 281; 581; 814; 274; 254; 854; 817; 547; 527; 827 & .3814 \\
451; 721; 218; 518; 184; 724; 245; 845; 187; 457; 572; 872 & .0992 \\
154; 127; 812; 815; 481; 427; 542; 548; 781; 754; 275; 278 & .0819 \\
415; 712; 128; 158; 148; 742; 425; 485; 178; 475; 752; 782 & .0194 \\
145; 172; 182; 185; 418; 472; 452; 458; 718; 745; 725; 728 & .0172 \\
\hline
214,241; 517,571; 251,521; 284,824; 847,874; 587,857 & .3965 \\
124,421; 157,751; 215,512; 248,842; 487,784; 578,875 & .0885 \\
142,412; 175,715; 125,152; 428,482; 478,748; 758,785 & .0149 \\
\hline
147,174,417,471,714,741; 258,285,528,582,825,852 & $1/6$ \\
\hline
\hline
\end{tabular}
\medskip
\caption{Three-way races modulo $q=9$}
\label{table3way9}
\end{table}

Again, the distribution of the primes into residue classes modulo 10 is
determined by their distribution mod 5, so the next modulus of interest
is $q=11$. Table~\ref{table2way11}
\begin{table}[b]
\begin{tabular}{||c|c||}
\hline 
\hline
${a_1a_2}$ & $\delta_{11; a_1,a_2}$ \\
\hline
\hline
23,25,64,69,71,75,81,89,T3,T4 & .761121 \\
\hline
21,24,61,63,73,79,84,85,T5,T9 & .731135 \\
\hline
29,65,74,83,T1 & .713943\\
\hline
NN,SS & $1/2$ \\
\hline
1T,38,47,56,92 & .286057 \\
\hline
12,16,36,37,42,48,58,5T,97,9T & .268865 \\
\hline
17,18,32,3T,46,4T,52,57,96,98 & .238879 \\
\hline
\hline
\end{tabular}
\medskip
\caption{Two-way races modulo $q=11$}
\label{table2way11}
\end{table}

shows the calculated densities for the two-way races modulo 11, where
we have used the symbol T to represent the residue 10 mod 11. In the
middle row, the entry ``NN,SS'' refers to the forty pairs $\{a_1,a_2\}$
where $a_1$ and $a_2$ are either both among the nonsquares
$\{2,6,7,8,{\rm T}\}$ or both among the squares $\{1,3,4,5,9\}$ mod 11.

We do not include the calculations of the three-way races mod 11 for
reasons of space. Using Theorem 2 it can be checked that of the 120
distinct (unordered) triples of residues mod 11, the twenty triples of the
form $\{ab^{-1},a,ab\}$ with $b$ a square mod 11 comprise two isomorphism
classes of race games of ten triples each; a race game in either of these
isomorphism classes has only two distinct densities, one taken by four
permutations of the triple and the other taken by the other two
permutations. The forty triples of the form $\{ab^{-1},a,ab\}$ with $b$ a
nonsquare mod 11 form four isomorphism classes with ten triples in each
class; a race game in one of these classes has three distinct densities,
each taken by a pair of permutations with the same middle element.
Finally, the remaining sixty triples form three isomorphism classes of
twenty race games each; a race game in one of these classes has all six
densities distinct. There are 34 densities that remain to be calculated
after these symmetries from Theorem 2 are taken into account, and the
calculations reveal that these 34 densities are indeed distinct.

As mentioned previously, determining the densities in a five-way race
game lies beyond the scope of the computing methods used for the
calculations in this paper (though this barrier is only technological,
as Theorem 4 is valid for arbitrarily large race games). If this
barrier were overcome
(for example, by recoding in a lower level computing language)
the five-way race among the squares mod 11 and
the five-way race among the nonsquares mod 11 would be natural and
interesting questions to consider, especially in light of the
nearly-cyclic behavior of the leaders in these five-way race games
reported by Bays and Hudson~\cite{BH}. Because of the symmetries of
Theorem 2, it turns out that only eight distinct densities would need
to be calculated for both of these five-way race games to be
completely determined.

Tables~\ref{table2way12},~\ref{table3way12}, and~\ref{table4way12}
\begin{table}[b]
\begin{tabular}{||c|c||}
\hline \hline
${a_1a_2}$ & $\delta_{12; a_1,a_2}$ \\ 
\hline
\hline
51 & .9992059 \\
15 & .0007941 \\
\hline
71 & .9986061 \\
17 & .0013939 \\
\hline
E1 & .9999766 \\
1E & .0000234 \\
\hline
57,5E,75, & \\
\;7E,E5,E7 & \up{$1/2$} \\
\hline
\hline
\end{tabular}
\medskip
\caption{Two-way races modulo $q=12$}
\label{table2way12}
\end{table}

show the two-way, three-way,
and four-way race games modulo 12, respectively, using the symbol E to
represent the residue 11 mod 12.
\begin{table}[t]
\begin{tabular}{||c|c||c|c||c|c||c|c||}
\hline \hline
$a_1a_2a_3$ & $\delta_{12; a_1,a_2,a_3}$ &
$a_1a_2a_3$ & $\delta_{12; a_1,a_2,a_3}$ &
$a_1a_2a_3$ & $\delta_{12; a_1,a_2,a_3}$ &
$a_1a_2a_3$ & $\delta_{12; a_1,a_2,a_3}$ \\
\hline
\hline
751 & $.4992728$ & 5E1 & $.4999772$ & 7E1 & $.4999780$ & & \\
571 & $.4986582$ & E51 & $.4992062$ & E71 & $.4986066$ & \up{57E,E75} & \up{$.1984521$} \\
\cline{7-8} 
517 & $.0012750$ & E15 & $.0007931$ & E17 & $.0013919$ & & \\
715 & $.0006751$ & 51E & $.0000225$ & 71E & $.0000214$ & \up{E57,75E} & \up{$.1799849$} \\
\cline{7-8}
157 & $.0000668$ & 1E5 & $.0000006$ & 1E7 & $.0000015$ & & \\
175 & $.0000521$ & 15E & $.0000003$ & 17E & $.0000006$ & \up{5E7,7E5} & \up{$.1215630$} \\
\hline
\hline 
\end{tabular}
\medskip
\caption{The four three-way races modulo $q=12$}
\label{table3way12}
\end{table}

Since the multiplicative group mod 12 is isomorphic to the multiplicative
group mod 8 (both groups being isomorphic to the Klein group of order 4), the
various symmetries present in Tables~\ref{table2way12},~\ref{table3way12},
and~\ref{table4way12} mirror those found in
Tables~\ref{table2way8},~\ref{table3way8}, and~\ref{table4way8}. As with the
modulus 8 case, all the characters mod 12 are real-valued, and so we feel
justified in reporting seven decimal places of the numbers in these
tables.
\begin{table}[t]
\begin{tabular}{||c|c|c|c|c|c|c|c||}
\hline \hline
$a_1a_2a_3a_4$ & $\delta_{12; a_1,a_2,a_3,a_4}$ &
$a_1a_2a_3a_4$ & $\delta_{12; a_1,a_2,a_3,a_4}$ &
$a_1a_2a_3a_4$ & $\delta_{12; a_1,a_2,a_3,a_4}$ &
$a_1a_2a_3a_4$ & $\delta_{12; a_1,a_2,a_3,a_4}$ \\ \hline \hline
157E & ${}<10^{-7}$ & 517E & $.0000004$ & 715E & $.0000001$ & E157 & $.0000664$ \\
15E7 & $.0000001$ & 51E7 & $.0000010$ & 71E5 & $.0000002$ & E175 & $.0000519$ \\
175E & ${}<10^{-7}$ & 571E & $.0000152$ & 751E & $.0000059$ & E517 & $.0011332$ \\
17E5 & ${}<10^{-7}$ & 57E1 & $.1984364$ & 75E1 & $.1799788$ & E571 & $.1787850$ \\
1E57 & $.0000002$ & 5E17 & $.0001403$ & 7E15 & $.0000243$ & E715 & $.0006505$ \\
1E75 & $.0000001$ & 5E71 & $.1214216$ & 7E51 & $.1215384$ & E751 & $.1977496$ \\ \hline \hline 
\end{tabular}
\medskip
\caption{The full four-way race modulo $q=12$}
\label{table4way12}
\end{table}

Notice from Table~\ref{table3way12} that the densities
$\delta_{12;5,11,1}$ and $\delta_{12;7,11,1}$ only differ by one unit
in the sixth decimal place, and that there are several other entries
that differ by similarly small amounts owing to their small
size. Nevertheless, we see no reason to believe that any of the
twenty-one densities in Table~\ref{table3way12} is equal to any
another. Similar remarks hold for the twenty-four densities in
Table~\ref{table4way12} and for the corresponding
Tables~\ref{table3way8} and~\ref{table4way8} for the race games modulo
8. One observation supporting our view is that whenever the symmetries
of Theorem 2 imply that two densities are equal, the computed
densities agree to within a few multiples of the default machine
precision rather than to only five or six decimal places.

\bigskip
\section{Proofs of Theorems 2 and 3.}
\setcounter{equation}{0}

In this section we establish Theorem 2, concerning symmetries of the
densities $\delta_\qaa$ under certain permutations of the residue
classes $\{a_1,\dots,a_r\}$, and Theorem 3, giving some strict
inequalities in the same setting. We first present the proof of
Theorem 3 since it is somewhat simpler than that of Theorem 2.

\medskip
\noindent{\bf Proof of Theorem 3.}
Let $a_1$, $a_2$, and $a_3$ be distinct reduced
residue classes mod $q$. We begin with the simple observation that if
$x$ is a real number such that $\pi(x;q,a_1)>\pi(x;q,a_2)$, then the
quantity $\pi(x;q,a_3)$ must either equal one of $\pi(x;q,a_1)$ and
$\pi(x;q,a_2)$, lie between them, exceed both, or be exceeded by
both. This observation leads to the density identity
\begin{equation}
\delta_{q;a_1,a_2} = \delta_{q;a_3,a_1,a_2} + \delta_{q;a_1,a_3,a_2}
+ \delta_{q;a_1,a_2,a_3},
\label{a3choices}
\end{equation}
since the set of real numbers $x$ such that
$\pi(x;q,a_3)=\pi(x;q,a_1)$ or $\pi(x;q,a_3)=\pi(x;q,a_2)$ has
density zero, as mentioned in Section 2.1. It follows that
\begin{equation}
\delta_{q;a_1,a_2,a_3} - \delta_{q;a_3,a_2,a_1} = \delta_{q;a_1,a_2} -
\delta_{q;a_3,a_2},
\label{twotermsleft}
\end{equation}
by using the appropriate identity of the type~(\ref{a3choices}) on
both terms on the right-hand side of
equation~(\ref{twotermsleft}) and simplifying.

Now we can use our knowledge of the two-way densities to study the
difference on the left-hand side of~(\ref{twotermsleft}). In
particular, if $c(q,a_1)=c(q,a_2)$ then $\delta_{q;a_1,a_2}=1/2$, and
hence $\delta_{q;a_1,a_2,a_3} - \delta_{q;a_3,a_2,a_1} = 1/2 -
\delta_{q;a_3,a_2}$, an expression whose sign is known
from the work of Rubinstein and Sarnak. More
specifically, if $N$ and $N'$ are nonsquares mod $q$ while $S$ is a
square mod $q$, then $\delta_{q;N,N',S} - \delta_{q;S,N',N} =
1/2 - \delta_{q;S,N'} > 0$;
therefore $\delta_{q;N,N',S}>\delta_{q;S,N',N}$, which establishes
part (a) of the theorem. Similarly, if $N$ is a nonsquare mod $q$
while $S$ and $S'$ are squares mod $q$, then
$\delta_{q;S',S,N}<\delta_{q;N,S,S'}$, which establishes part (b) of
the theorem.

Another application is to the difference
$\delta_{q;N,S,N'}-\delta_{q;N',S,N}$ when $N$ and $N'$ are
nonsquares mod $q$ while $S$ is a square mod $q$. In this case
equation~(\ref{twotermsleft}) becomes
\begin{equation*}
\delta_{q;N,S,N'} - \delta_{q;N',S,N} = \delta_{q;N,S} -
\delta_{q;N',S},
\end{equation*}
which immediately implies part (c) of the theorem. The analogous
observation about the difference
$\delta_{q;S,N,S'}-\delta_{q;S',N,S}$ when
$S$ and $S'$ are squares mod $q$ while $N$ is a nonsquare mod $q$
establishes part (d) of the theorem.\qed\medskip

We remark that the identity~(\ref{twotermsleft}), applied when $a_1$,
$a_2$, and $a_3$ are all nonsquares mod $q$, becomes
$\delta_{q;a_1,a_2,a_3} - \delta_{q;a_3,a_2,a_1} = 0$; this is
another way of seeing that the densities calculated in Theorem 1 are
equal in pairs as indicated.

Our next goal is to establish Theorem 2. Before doing so it will be
helpful to recall the relationships between the density $\delta_\qaa$
and the measures $\mu_\qaa$ and $\rho_\qaa$. We begin by recalling
from equation~(\ref{deltaasmu}) that
\begin{equation}
\delta_\qaa = \mathop{\int \dots \int}_{x_1 > \dots > x_r} d \mu_\qaa
\, .
\label{deltaasmureprise}
\end{equation}
We remark that if $\sigma$ is a permutation of the indices
$\{1,\dots,r\}$, then we can express the density
$\delta_{q;a_{\sigma(1)},\dots,a_{\sigma(r)}}$ in two different ways:
we have
\begin{equation*}
\delta_{q;a_{\sigma(1)},\dots,a_{\sigma(r)}} = \mathop{\int \dots
\int}_{x_1 > \dots > x_r} d \mu_{q;a_{\sigma(1)},\dots,a_{\sigma(r)}}
\end{equation*}
corresponding to the formula~(\ref{deltaasmureprise}), but we also
have the alternate form
\begin{equation*}
\delta_{q;a_{\sigma(1)},\dots,a_{\sigma(r)}} = \mathop{\int \dots
\int}_{x_{\sigma(1)} > \dots > x_{\sigma(r)}} d \mu_\qaa
\end{equation*}
since $\mu_\qaa$ is the limiting distribution of the vector
$(E(x;q,a_1),\dots,E(x;q,a_r))$, whose coordinated are ordered by size
exactly as the coordinates of the vector
$(\pi(x;q,a_1),\dots,\pi(x;q,a_r))$.

If we make the change of variables $u_1 = x_1 - x_2$, \dots,
$u_{r-1}=x_{r-1} - x_r$, $u_r=x_r$ and integrate out the variable
$u_r$, as in Section 2.5, the formula~(\ref{deltaasmureprise}) becomes
\begin{equation}
\delta_\qaa = \mathop{\int \dots \int}_{u_1>0,\, \dots,\, u_{r-1}>0} d
\rho_\qaa.
\label{byvirtueof}
\end{equation}
For the special permutation $\sigma$ that reverses the set
$\{1,\dots,n\}$, we see that
\begin{equation*}
x_{\sigma(1)} > \dots > x_{\sigma(r)} \quad\Longleftrightarrow\quad
x_r > \dots > x_1 \quad\Longleftrightarrow\quad
u_{r-1}<0,\, \dots,\, u_1<0.
\end{equation*}
Consequently, we have
\begin{equation}
\delta_{q;a_r,\dots,a_1} = \mathop{\int \dots \int}_{u_1<0,\, \dots,\,
u_{r-1}<0} d \rho_\qaa
\label{byvirtuereversed}
\end{equation}
as a companion formula to equation~(\ref{byvirtueof}).

As a final prerequisite to the proof of Theorem 2 we recall from
equation~(\ref{muhatformula}) the explicit formula
\begin{equation}
\hat\mu_\qaa(\xi_1,\dots.\xi_r)
= \exp \bigg( i \sum_{j=1}^r c(q,a_j)\xi_j \bigg) \prod
\begin{Sb}\chi\bmod q \\ \chi\ne\chi_0\end{Sb} F\bigg( \bigg|
\sum_{j=1}^r \chi(a_j)\xi_j \bigg|, \chi \bigg),
\label{muhatreprise}
\end{equation}
for the Fourier transform of $\mu_\qaa$, and the related
formula~(\ref{rhohatformula})
\begin{multline}
\hat\rho_\qaa(\eta_1,\dots,\eta_{r-1}) = \exp\bigg( \sum_{j=1}^{r-1}
(c(q,a_j)-c(q,a_{j+1}))\eta_j \bigg) \\
\times \prod \begin{Sb}\chi\bmod q \\ \chi\ne\chi_0\end{Sb} F\bigg(
\bigg| \sum_{j=1}^{r-1} (\chi(a_j)-\chi(a_{j-1}))\eta_j \bigg|, \chi
\bigg)
\label{keyfact}
\end{multline}
for the Fourier transform of $\rho_\qaa$.

\medskip
\noindent{\bf Proof of Theorem 2.}
Let $a_j^{-1}$ denote the inverse of $a_j$ mod $q$. We will show that
the Fourier transforms $\hat\mu_\qaa$ and
$\hat\mu_{q;a_1^{-1},\dots,a_r^{-1}}$ are the same function. This is
enough to establish part (a), since the densities $\mu_\qaa$ and
$\mu_{q;a_1^{-1},\dots,a_r^{-1}}$ will then be identical, which by
equation~(\ref{deltaasmureprise}) will imply
\begin{equation*}
\delta_\qaa = \mathop{\int \dots \int}_{x_1 > x_2 > \dots > x_r}
d\mu_\qaa = \mathop{\int \dots \int}_{x_1 > x_2 > \dots > x_r}
d\mu_{q;a_1^{-1},\dots,a_r^{-1}} = \delta_{q;a_1^{-1},\dots,a_r^{-1}}.
\end{equation*}

We use the formula~(\ref{muhatreprise}) for $\hat\mu_\qaa$ and the
analogous formula for $\hat\mu_{q;a_1^{-1},\dots,a_j^{-1}}$. Notice
that the square roots of $a_j^{-1}$ are precisely the inverses mod
$q$ of the square roots of $a_j$. In particular,
$c(q,a_j^{-1})=c(q,a_j)$, and so the exponential term in the
formula~(\ref{muhatreprise}) is unchanged if we replace each $a_j$ by
$a_j^{-1}$. Moreover, we see that for each character $\chi$ mod $q$,
\begin{equation*}
\bigg| \sum_{j=1}^r \chi(a_j^{-1})\xi_j \bigg| = 
\bigg| \sum_{j=1}^r \overline{\chi(a_j)} \xi_j \bigg| = 
\bigg| \sum_{j=1}^r \chi(a_j)\xi_j \bigg|
\end{equation*}
since the $\xi_j$ are real, so that each term $F(\cdot,\chi)$
in~(\ref{muhatreprise}) is also unchanged by replacing all of the
$a_j$ with the $a_j^{-1}$. This shows that $\hat\mu_\qaa =
\hat\mu_{q;a_1^{-1},\dots,a_r^{-1}}$, which establishes part (a) of
the theorem.

We use a similar strategy to prove part (b). Let $b$ be a reduced
residue class mod $q$ such that $c(q,a_j)=c(q,ba_j)$ for each $1\le
j\le r$. Because of this hypothesis, the exponential term in the
formula~(\ref{muhatreprise}) is unchanged if we replace each $a_j$ by
$ba_j$ as above. Moreover, for each character $\chi$ mod~$q$,
\begin{equation}
\bigg| \sum_{j=1}^r \chi(ba_j)\xi_j \bigg| =
\bigg| \chi(b) \sum_{j=1}^r \chi(a_j)\xi_j \bigg| =
|\chi(b)| \bigg| \sum_{j=1}^r \chi(a_j)\xi_j \bigg| =
\bigg| \sum_{j=1}^r \chi(a_j)\xi_j \bigg|,
\label{partsbandc}
\end{equation}
so that each term $F(\cdot,\chi)$ in~(\ref{muhatreprise}) is also
unchanged by replacing all of the $a_j$ with the $ba_j$. This
shows that $\hat\mu_\qaa = \hat\mu_{q;ba_1,\dots,ba_r}$, which
establishes part (b) of the theorem.

The proofs of parts (c) and (d) rely on the formula~(\ref{keyfact})
for the function $\hat\rho_\qaa$. When the $a_j$ are all squares mod
$q$, then the exponential term in~(\ref{keyfact}) is identically
1. Moreover, if $b$ is a square mod $q$ then each $ba_j$ is also a
square, while if $b$ is a nonsquare mod $q$ then each $ba_j$ is a
nonsquare; in either case we have $c(q,ba_1)=\dots=c(q,ba_r)$, so that
the exponential term in the analogous formula to
equation~(\ref{keyfact}) for $\hat\rho_{q;ba_1,\dots,ba_r}$ is also
identically 1. Since the chain of equalities~(\ref{partsbandc}) again
shows that each term $F(\cdot,\chi)$ is unchanged upon replacing the
$a_j$ with $ba_j$, we see that
$\hat\rho_\qaa=\hat\rho_{q;ba_1,\dots,ba_r}$ and so
$\delta_\qaa=\delta_{q;ba_1,\dots,ba_r}$ by virtue of
equation~(\ref{byvirtueof}), which establishes part (c) of the
theorem.

For part (d) we begin with the formula~(\ref{byvirtuereversed}) for
$\delta_{q;a_r,\dots,a_1}$. As remarked above, the exponential term of
$\hat\rho_\qaa$ is identically 1 when the $a_j$ are all squares mod
$q$, so that $\hat\rho_\qaa$ will be real valued. Since $\rho_\qaa$
is real-valued as well, we conclude that $\rho_\qaa$ is symmetric
through the origin. Hence making the change of variables $u_j \mapsto
-u_j$ for each $1\le j\le r$ in equation~(\ref{byvirtuereversed}), we
obtain
\begin{equation*}
\delta_{q;a_r, \dots , a_1} = \mathop{\int \dots \int}_{u_1>0,
\dots , u_{r-1}>0} d \rho_\qaa = \delta_\qaa,
\end{equation*}
which establishes part (d) of the theorem.

To establish part (e), we first consider the relationship between
$\hat\rho_\qaa$ and $\hat\rho_{q;ba_1,\dots,ba_r}$ (note that the
residue classes $ba_j$ have not yet been reversed in the second
subscript). Again, equation~(\ref{partsbandc}) shows that replacing
each $a_j$ with $ba_j$ does not change the terms of the form
$F(\cdot,\chi)$, and so we only need to consider the exponential
term. Because the quantity $c(q,a)$ can only take the two values
${-1}$ and $c(q,1)$, we see that if $c(q,a')\ne c(q,a)$ then
$c(q,a')=c(q,1)-1-c(q,a)$. It follows that under our hypothesis that
$c(q,ba_j)\ne c(q,a_j)$ for each $1\le j\le r$; but we also have
\begin{equation*}
c(q,ba_{j+1})-c(q,ba_j)=-(c(q,a_{j+1})-c(q,a_j)),
\end{equation*}
and so the imaginary expression in the exponential term in
equation~(\ref{muhatreprise}) is negated upon replacing each $a_j$ by
$ba_j$. The end result is that
$\hat\rho_{q;ba_1,\dots,ba_r}=\overline{\hat\rho_\qaa}$, which
implies that when the measure $\rho_{q;ba_1,\dots,ba_r}$ is reflected
through the origin, the resulting measure is identical to
$\rho_\qaa$.

Since we can express
\begin{equation*}
\delta_{q;ba_r,\dots,ba_1} =
\mathop{\int\dots\int}_{u_1<0,\,\dots,\,u_{r-1}<0}
\rho_{q;ba_1,\dots,ba_r}
\end{equation*}
as in equation~(\ref{byvirtuereversed}), we can make the change of
variables $u_j\mapsto-u_j$ for each $1\le j\le r-1$ to see that
\begin{equation*}
\delta_{q;ba_r,\dots,ba_1} =
\mathop{\int\dots\int}_{u_1>0,\,\dots,\,u_{r-1}>0}
\rho_\qaa = \delta_\qaa.
\end{equation*}
This establishes the final assertion of the theorem.\qed

\bigskip
\section{Remarks, Questions, and Open Problems.}
\setcounter{equation}{0}

In this final section, we collect together several observations,
unanswered questions, and conjectures concerning the results of this paper.

\medskip\noindent{\it Systems of inequalities with one equality.\/}
Since we know that $\delta_{q;a,b}$ and $\delta_{q;b,a}$ are both positive
(assuming GRH and LI), it follows that each inequality
$\pi(x;q,a)>\pi(x;q,b)$ and $\pi(x;q,b)>\pi(x;q,a)$ has arbitrarily
large solutions, and therefore $\pi(x;q,a)=\pi(x;q,b)$ for
infinitely many integers $x$. However, knowing that $\delta_{q;a,b,c}$
and $\delta_{q;b,a,c}$ are both positive---i.e., that each string of
inequalities
\begin{equation*}
\pi(x;q,a)>\pi(x;q,b)>\pi(x;q,c) \quad\hbox{and}\quad
\pi(x;q,b)>\pi(x;q,a)>\pi(x;q,c)
\end{equation*}
has arbitrarily large solutions---does not imply that there are
necessarily any solutions to
$\pi(x;q,a)=\pi(x;q,b)>\pi(x;q,c)$. Undoubtably, the equality
$\pi(x;q,a)=\pi(x;q,b)$ should hold infinitely often both when their
common value exceeds $\pi(x;q,c)$ and when their value is exceeded by
$\pi(x;q,c)$. We conjecture more generally that for any given integer
$1\le j\le r$ and reduced residue classes $a_1$, \dots, $a_r$ and
$a_j'$ mod $q$, the conditions

\medskip
\centerline{
\begin{tabular}{ccc}
$\pi(x;q,a_1)>\dots>$ & \hskip-.5em$\pi(x;q,a_j)$ &
\hskip-.5em$>\dots>\pi(x;q,a_r)$ \\
& \hskip-.5em$||$ & \\
& \hskip-.5em$\pi(x;q,a'_j)$ \\
\end{tabular}
}
\medskip

\noindent should be satisfied for infinitely many integers $x$.

\medskip\noindent{\it Multiple equalities.\/}
Another direction along these lines involves solutions to
\begin{equation}
\pi(x;q,a_1)=\pi(x;q,a_2)=\dots=\pi(x;q,a_r)  \label{allpisequal}
\end{equation}
when $r\ge3$. If we consider the vectors
\begin{equation}
V_\qaa(n) = \big ( \pi(p_n;q,a_1)-\pi(p_n;q,a_2 ) , \dots,
\pi(p_n;q,a_{r-1})-\pi(p_n;q,a_r) \big ) ,
\label{vvectors}
\end{equation}
where $p_n$ denotes the $n$th prime, then the sequence of vectors
$\{V_\qaa(n)\}$ might reasonably be expected to resemble a random walk on
$\Z^{r-1}$, where the possible steps at each stage are $(1,0,\dots,0)$,
$(-1,1,0,\dots,0)$, \dots, $(0,\dots,0,-1,1)$, and $(0,\dots,0,-1)$ and are
chosen with roughly equal probabilities. (Even though the Chebyshev bias will
cause a drift in the mean behavior of the vectors~(\ref{vvectors}), this
drift has the same order of magnitude as the standard deviation of the
random walk).

Since random walks on $\Z^n$ return to any point infinitely often with
probability 1 when $n=1$ or 2 but fail to do so with probability 1
when $n\ge3$ (Polya~\cite{polya}), this heuristic leads to the
prediction that the system of equalities~(\ref{allpisequal}) has
infinitely many solutions when $r\le3$ but only finitely many
solutions for $r\ge4$.  Similar reasoning suggests that any pair of
equalities
\begin{equation*}
\pi(x;q,a_1)=\pi(x;q,a_2), \quad \pi(x;q,a_3)=\pi(x;q,a_4)
\end{equation*}
with $a_1$, \dots, $a_4$ distinct should simultaneously hold for
arbitrarily large values of $x$, but three or more equalities will
hold simultaneously only finitely many times.  Further, we might
expect that the conditions

\medskip
\centerline{
\begin{tabular}{ccc}
$\pi(x;q,a_1)>\dots>$ & \hskip-.5em$\pi(x;q,a_j)$ &
\hskip-.5em$>\dots>\pi(x;q,a_r)$ \\
& \hskip-.5em$||$ & \\
& \hskip-.5em$\pi(x;q,a'_j)$ \\
& \hskip-.5em$||$ & \\
& \hskip-.5em$\pi(x;q,a''_j)$ \\
\end{tabular}
}
\medskip

\noindent and

\medskip
\centerline{
\begin{tabular}{ccccc}
$\pi(x;q,a_1)>\dots>$ & \hskip-.5em$\pi(x;q,a_i)$ &
\hskip-.5em$>\dots>$ & \hskip-.5em$\pi(x;q,a_j)$ &
\hskip-.5em$>\dots>\pi(x;q,a_r)$ \\
& \hskip-.5em$||$ & & \hskip-.5em$||$ & \\
& \hskip-.5em$\pi(x;q,a'_i)$ & & \hskip-.5em$\pi(x;q,a'_j)$ \\
\end{tabular}
}
\medskip

\noindent should hold for infinitely many integers $x$, but that analogous
conditions involving three or more equalities would not.

\medskip\noindent{\it Bias factors.\/}
To try to quantify the Chebyshev biases for $r$-tuples of reduced residue
classes $a_j\bmod q$ for all $r\ge2$, let us define the ``bias factor''
$\beta_\qaa$ to be the difference between the number of nonsquares
preceding squares among the $a_j$ and the number of squares preceding
nonsquares:
\begin{equation}
\begin{split}
\beta_\qaa &= \#\{ i<j\colon a_i\ne {\hskip-1em\qed},\, a_j=
{\hskip-1em\qed} \} - \#\{ i<j\colon a_i= {\hskip-1em\qed},\, a_j\ne
{\hskip-1em\qed} \} \\
&= \mathop{\sum\sum}_{1\le i<j\le r} {c(q,a_j)-c(q,a_i) \over
c(q,1)+1} \\
&= {1\over c(q,1)+1} \sum_{1\le j\le r} (2j-r-1)c(q,a_j).
\end{split}
\label{biasfactordef}
\end{equation}
For instance, when $r=2$ the possible bias factors are $\beta_{q;N,S}=1$,
$\beta_{q;N,N}=\beta_{q;S,S}=0$, and $\beta_{q;S,N}=-1$. Rubinstein and
Sarnak proved that the sign of $\delta_{q;a,b}-1/2$ equals the sign of
$\beta_{q;a,b}$ in this notation, thereby showing that
\begin{equation*}
\beta_{q;a,b} > \beta_{q;a',b'} \implies \delta_{q;a,b} > \delta_{q;a',b'}.
\end{equation*}
The converse to this statement is false: the first two lines of
Table~\ref{table2way7} shows that $\delta_{q;a,b}$ and $\delta_{q;a',b'}$
can be different even when $\beta_{q;a,b}=\beta_{q;a',b'}$, for instance.

We might hope that the bias factors $\beta_\qaa$ would provide some
information about the relative sizes of the $\delta_\qaa$, perhaps in the form
of the implication
\begin{equation}
\beta_\qaa > \beta_{q;b_1,\dots,b_r} \implies \delta_\qaa >
\delta_{q;b_1,\dots,b_r}
\label{biasfactorconj}
\end{equation}
for any fixed $r$. In this regard, it is worth remarking that all of
the symmetries in Theorem 2 are equalities between two $r$-tuples of
residues with equal bias factors.  Examining the densities computed in
Section 4, we observe that the implication~(\ref{biasfactorconj})
holds most of the time, but we do note the following two anomalies:
\begin{itemize}
\item $\beta_{8;5,1,3,7}=\beta_{8;5,1,7,3}=-1>-3=\beta_{8;1,3,7,5}$, but
it appears from Table~\ref{table4way8} that $\delta_{8;1,3,7,5}$ slightly exceeds
both $\delta_{8;5,1,3,7}$ and $\delta_{8;5,1,7,3}$;
\item $\beta_{12;7,1,11,5}=\beta_{12;7,1,5,11}=-1>-3=\beta_{12;1,11,5,7}$,
but it appears from Table~\ref{table4way12} that $\delta_{12;1,11,5,7}$
slightly exceeds both $\delta_{12;7,1,11,5}$ and $\delta_{12;7,1,5,11}$.
\end{itemize}
It would therefore be of interest, in connection with determining whether
the implication~(\ref{biasfactorconj}) is always valid, to compute more
precisely the densities just mentioned in order to verify the apparent
inequalities.

Unfortunately, the computation of the densities to arbitrary precision is
not simply a matter of reducing $\epsilon$ and increasing $C$ and letting
a bigger computer run for a longer period of time. The major source of error in
these computations is the effect of truncating the infinite product defining
the functions $F(z,\chi)$ to form the approximations $F_T(z,\chi)$ (see
Section 3.3); to decrease this error it would be necessary to compute zeros
of the relevant $L$-functions to a height greater than 10,000, and perhaps
to greater precision than twelve decimal places as well.

It is certainly conceivable that some definition of bias factor
different from~(\ref{biasfactordef}) might be better suited to the
role of $\beta_\qaa$, although it is hard to imagine what natural
definition would be able to explain the apparent anomalies noted
above. It might also be the case that the
implication~(\ref{biasfactorconj}) is valid in more limited
settings---for instance, when we restrict to $r$-tuples
$\{a_1,\dots,a_r\}$ and $\{b_1,\dots,b_r\}$ where exactly half of the
$a_j$ are nonsquares and half squares, and similarly for the $b_j$.

\medskip\noindent{\it Convergence to unbiased distribution.\/}
Rubinstein and Sarnak proved \cite[Theorem 1.5]{RS} that for a fixed
integer $r\ge2$,
\begin{equation}
\big( \max_{a_1,\dots,a_r} |r!\delta_\qaa - 1| \big) \to 0
\label{maxtozero}
\end{equation}
as $q$ tends to infinity (where the maximum is taken over all $r$-tuples of
distinct reduced residue classes mod $q$), so that biases of any sort become
less and less evident with increasing moduli. Thus although the biases in the
two-way races mod 8 and mod 12 are more pronounced than those in the two-way
races mod 4, 5, and 7 owing to the larger values of $c(8,1)=c(12,1)=3$, these sorts
of extreme biases will not continue (even with a sequence of
moduli such as $q_n = 4p_2p_3\dots p_n$, say, which satisfies
$c(q_n,1)=2^n-1$).

On the other hand, it might happen that an extremely negatively biased
density such as $\delta_{q;S_1,\dots,S_n,N_1,\dots,N_n}$ might tend to
zero much more rapidly than $1/(2n)!$\ as $n$ increases, while an
extremely positively biased density such as
$\delta_{q;N_1,\dots,N_n,S_1,\dots,S_n}$ might behave more like
$1/(n!)^2$. In general, one could investigate the uniformity of the
statement~(\ref{maxtozero}), i.e., attempt to show that the statement
holds uniformly for all $r\le r_0$ for some integer-valued function
$r_0=r_0(q)$ satisfying $2\le r_0\le\phi(q)$. For instance, is it the
case that
\begin{equation}
\limsup_{q\to\infty} \big( \max_{a_1,\dots,a_{r_0}}
r_0!\,\delta_{q;a_1,\dots,a_{r_0}} \big) = \infty \quad\hbox{and}\quad
\liminf_{q\to\infty} \big( \min_{a_1,\dots,a_{r_0}}
r_0!\,\delta_{q;a_1,\dots,a_{r_0}} \big) = 0
\end{equation}
if $r_0 = r_0(q)$ grows sufficiently quickly, and if so, how quickly
must $r_0$ grow with $q$ for these phenomena to emerge? We certainly
conjecture that
\begin{equation}
\lim_{q\to\infty} \big( \max_{a_1,\dots,a_{r_0}}
r_0!\,\delta_{q;a_1,\dots,a_{r_0}} \big) = 0
\end{equation}
for any arbitrary function $r_0=r_0(q)$ tending to infinity with $q$,
but at this point it seems nontrivial to prove this modest result even
for $r_0=\phi(q)$ itself.
\renewcommand{\labelenumi}{(\arabic{enumi})}

\medskip\noindent{\it Race-game symmetries, isomorphisms, and order
equivalences.\/}
Another question of interest is whether there exist more symmetry
results of the type arising in Theorem 2. Reviewing the proof of
Theorem 2, we see that all of the symmetries therein are consequences
of provable equalities between two distributions of the type
$\mu_\qaa$ or $\rho_\qaa$ (possibly after reflecting one of the
distributions through the origin). We can then ask
\begin{enumerate}
\item whether there exist any equalities between these distributions other
than those used in the proof of Theorem 2;
\item whether there can be numerical ``coincidences'' between two densities
even though their underlying distributions are not related.
\end{enumerate}
An answer to question (1) might be forthcoming from a careful analysis of the
Fourier transforms $\hat\rho_\qaa$ of the distributions $\rho_\qaa$.
As for question (2), it
seems reasonable to believe the phenomenon addressed therein can never
occur, but proving such a claim seems very difficult.

In support of the possibility that Theorem 2 accounts for all numerical
equalities between the densities $\delta_\qaa$, we remark that among the
densities computed in Section 4, each time a symmetry from Theorem 2
was applicable the corresponding computed densities were equal to within a small multiple of
the machine precision. Conversely, all such numerical equalities observed
among the computed densities are accounted for by the symmetries already
asserted in Theorem~2.

Symmetries among individual densities $\delta_\qaa$ are of course closely
related to isomorphisms between complete race games of $r$-tuples. Theorem
2 implies that the following bijections between $r$-tuples induce
isomorphisms of race games:
\begin{itemize}
\item the map $\tau(a_j)\equiv a_j^{-1}\pmod q$ between the $r$-tuples 
$\{a_1,\dots,a_r\}$ and $\{a_1^{-1},\dots,a_r^{-1}\}$;
\item the map $\tau(a_j)\equiv ba_j\pmod q$ between the $r$-tuples 
$\{a_1,\dots,a_r\}$ and $\{ba_1,\dots,ba_r\}$, if either
$c(q,a_j)=c(q,1)$ for each $1\le j\le r$ or $c(q,a_j)=c(q,ba_j)$ for each
$1\le j\le r$;
\item the map $\tau(a_j)\equiv ba_{r+1-j}\pmod q$ between the $r$-tuples 
$\{a_1,\dots,a_r\}$ and $\{ba_r,\dots,ba_1\}$, if $c(q,a_j)\ne c(q,ba_j)$
for each $1\le j\le r$;
\item either bijection $\tau\colon \{a,b\}\to\{a',b'\}$, if $c(q,a)=c(q,b)$
and $c(q',a')=c(q',b')$;
\item any bijection $\tau\colon \{a,b,c\}\to\{a',b',c'\}$, if there exists
$\rho\not\equiv1\pmod q$ with $\rho^3\equiv1$ (mod~$q$) such that $b\equiv
a\rho\pmod q$ and $c\equiv a\rho^2\pmod q$ and an analogous
$\rho'\pmod{q'}$.
\end{itemize}
(Our definition of isomorphic race games required that the $r$-tuples
consist of reduced residues to the same modulus, but the definition has an
obvious extension to two $r$-tuples of residues to different moduli which
encompasses the last two isomorphisms.) We conjecture that any
isomorphism between two race games is induced by a composition of
bijections from this list; in particular, the only isomorphisms between race
games of distinct moduli are those race games with complete internal
symmetry, which were determined by Rubinstein and Sarnak.

A weaker relationship than isomorphic race games is {\it order-equivalent
race games\/}, where there exists a bijection $\tau$ on the set
$\{1,\dots,n\}$ such that
\begin{equation}
\delta_{q;a_{\sigma(1)},\dots,a_{\sigma(r)}} >
\delta_{q;a_{\sigma'(1)},\dots,a_{\sigma'(r)}} \quad\Longleftrightarrow\quad
\delta_{q';b_{\tau(\sigma(1))},\dots,b_{\tau(\sigma(r))}} >
\delta_{q';b_{\tau(\sigma'(1))},\dots,b_{\tau(\sigma'(r))}}
\label{orderequivdef}
\end{equation}
for any two permutations $\sigma$, $\sigma'$ of
$\{1,\dots,n\}$. Order-equivalent race games seem common for small
values of $r$. For instance, any two race games both of the form
$\{N,S\}$ are order-equivalent by Rubinstein and Sarnak's results. The
tables in Section 4 indicate many three-way race games that are
order-equivalent. The triples $\{N,N',1\}$ mod 7 with
$NN'\not\equiv-1$ mod 7, the triples $\{N,N',1\}$ mod 8, the triples
$\{N,N',1\}$ mod 9 with $NN'\not\equiv-1$ mod 9, and the triples
$\{N,N',1\}$ mod 12 are all order-equivalent to one another. Also, the
triples $\{N,-N^{-1},S\}$ mod 5, the triples $\{N,-N^{-1},S\}$ mod 7,
and the triples $\{N,-N^{-1},S\}$ mod~9 are all order-equivalent as
well (but note that these are not order-equivalent to the triples
$\{N,N,N\}$ mod 8 and mod 12).

We remark that from the values in Tables~\ref{table4way8}
and~\ref{table4way12}, the bijection
\begin{equation*}
\tau(1)=1, \quad\tau(3)=11, \quad\tau(5)=7, \quad\tau(7)=5
\end{equation*}
is quite close to inducing an order-equivalence between the full
four-way race games modulo 8 and 12, respectively (in the sense that the
values in these tables would only have to be modified by at most
$6\times10^{-5}$ in order for the condition~(\ref{orderequivdef}) to
always hold). It would certainly be interesting to try to establish
(or even classify) order-equivalent race games, especially for larger
values of $r$ and between $r$-tuples to different moduli.

\medskip\noindent{\it Another problem of Knapowski--Tur\'an.\/}
In their paper~\cite{KT}, Knapowski and Tur\'an pose many problems in
comparative prime number theory, several of which have been answered
by Rubinstein and Sarnak~\cite{RS} and herein.  We conclude by
mentioning one other problem given by Knapowski and Tur\'an
in~\cite{KT}. They ask whether, for any $r$-tuple $a_1,\dots,a_r$ of
reduced residue classes mod $q$, the inequalities
\begin{equation}
\pi(x;q,a_1) < {\li(x)\over\phi(q)},\quad \pi(x;q,a_2) <
{\li(x)\over\phi(q)},\quad \dots,\quad \pi(x;q,a_r) < {\li(x)\over\phi(q)}
\label{sixomega}
\end{equation}
simultaneously hold for arbitrarily large values of $x$. Each individual
inequality is unbiased if $a_j$ is a nonsquare mod $q$ and biased
negatively if $a_j$ is a square mod $q$. We remark here that if we
apply the method of Rubinstein and Sarnak~\cite{RS} to the error
term
\begin{equation*}
E_1(x;q,a) = {\log x\over\sqrt x} \big( \phi(q)\pi(x;q,a)-\li(x) \big),
\end{equation*}
which has been centered in a slightly different way than in the
definition~(\ref{Exqadef}) of $E(x;q,a)$, we can see that this
question of Knapowski and Tur\'an is answered in the affirmative,
and in fact the set of real numbers $x$ satisfying the
inequalities~(\ref{sixomega}) has positive density as well.

\bigskip
{\smaller\narrower

\centerline{\sc Acknowledgements}\medskip

The authors would like to gratefully acknowledge Robert Rumely, who
supplied us with the results of his calculations of the zeros of the
various $L(s,\chi)$ used in this paper, and also Kenneth Williams for
pointing out the closed-form expressions for $L'(1,\chi)$ cited
herein.  This work was supported by grants from the Natural Sciences
and Engineering Research Council of Canada.

}

\end{document}